\documentclass[twoside, 12pt]{paper}

\usepackage{a4wide, amsthm, amssymb, amscd, mathptmx}

\DeclareSymbolFont{rsfs}{OMS}{rsfs}{m}{n}
\DeclareSymbolFontAlphabet{\mathscr}{rsfs}
\DeclareSymbolFont{bbold}{U}{dsrom}{m}{n}
\DeclareSymbolFontAlphabet{\mathbb}{bbold}

\renewcommand{\Bbb}{\mathbb}
\renewcommand{\frak}{\mathfrak}
\renewcommand{\cal}{\mathscr}

\newcommand{\C}{{\Bbb C}}

\newcommand{\E}{{\cal E}}

\newcommand{\EL}{{\cal L}}
\newcommand{\End}{\mathop{\rm End}}

\newcommand{\F}{{\cal F}}
\newcommand{\G}{{\cal G}}

\newcommand{\GL}{\mathop{\rm GL}}

\newcommand{\Hom}{\mathop{\rm Hom}}

\newcommand{\id}{\mathop{\rm id}}

\newcommand{\n}{{\cal N}}
\newcommand{\Oh}{{\cal O}}

\renewcommand{\P}{{\Bbb P}}

\newcommand{\Q}{{\Bbb Q}}

\newcommand{\SL}{\mathop{\rm SL}}

\renewcommand{\theta}{\vartheta}
\newcommand{\Z}{{\Bbb Z}}
\newcommand{\N}{{\cal N}}
\newcommand{\la}{\lambda}
\newcommand{\lra}{\longrightarrow}

\newcommand{\lma}{\longmapsto}
\newcommand{\p}{\prime}
\newcommand{\q}{\quad}

\renewcommand{\phi}{\varphi}
\newcommand{\rk}{\mathop{\rm rk}}
\newcommand{\eps}{\varepsilon}

\newcommand{\ul}{\underline}
\newcommand{\ol}{\overline}

\theoremstyle{plain}
\newtheorem{Thm}{Theorem}[subsection]

\newtheorem{Prop}[Thm]{Proposition}
\newtheorem{Lem}[Thm]{Lemma}

\newtheorem*{Thm*}{Theorem}

\theoremstyle{remark}
\newtheorem{Rem}[Thm]{Remark}
\newtheorem*{Rem*}{Remark}

\newtheorem*{Conv}{Convention}
\newtheorem{Ex}[Thm]{\it Example}

\theoremstyle{definition}

\title{
Moduli for decorated tuples of sheaves and representation spaces for quivers\footnote{To appear in Proc.\ Indian
Acad.\ Sci.\ (Math Sci.).}}
\author{Alexander Schmitt
\institution{Universit\"at Duisburg-Essen\\ Campus Essen\\
Fachbereich 6: Mathematik und Informatik\\ D-45117 Essen\\ Deutschland
\\
\tt alexander.schmitt@uni-essen.de}}
\begin{document}
\maketitle
\begin{abstract}
We extend the scope of a former paper to vector bundle problems involving more than
one vector bundle. As the main application, we obtain the solution of the well-known
moduli problems of vector bundles associated with general quivers.
\end{abstract}
\section*{Introduction}
If we are given a projective manifold $X$, a reductive linear
algebraic group $G$, and a representation $\rho\colon G\lra
\GL(V)$, we may associate to every  principal $G$-bundle ${\cal
P}$ over $X$ a vector bundle ${\cal P}(V)$ with fibre $V$. The
objects of interest are pairs $({\cal P},\tau)$ where ${\cal P}$
is an algebraic principal $G$-bundle and $\tau\colon X\lra {\cal
P}(V)$ is a section of the associated vector bundle. Motivated,
e.g., by the quest for differentiable invariants of
$4$-manifolds, one associates to the data of $G$, $\rho$, and a
fixed differentiable principal $G$-bundle $P$ certain vortex
equations. Via a so-called Kobayashi-Hitchin correspondence, the
solutions of these vortex equations have an interpretation as
pairs $({\cal P},\tau)$ as above, satisfying certain stability
conditions which may be understood in purely algebraic terms.
Here, ${\cal P}$ is an algebraic structure on the bundle $P$. The
Kobayashi-Hitchin correspondence first arose in the context of
vector bundles, i.e., when no representation is given (see
\cite{LT}), and was then considered in various special cases
before Banfield \cite{Ba} gave it a unified treatise. It was
afterwards widely extended to more general contexts (\cite{MiR},
\cite{AG}, \cite{OT}, \cite{BGPMiR}, \cite{LT2}). In order to apply the
machinery of Algebraic Geometry to the gauge theoretic moduli
space for the pairs $({\cal P},\tau)$ with ${\cal P}$ of
topological type $P$ satisfying the stability conditions, one
must equip it with an algebraic structure and find a suitable
compactification. One is therefore led to a purely
algebro-geometric moduli problem. Another motivation to study
this kind of moduli problems that comes from within Algebraic
Geometry is the fact that many interesting classification
problems for projective manifolds may be encoded by data of the
above type. We will give an example below and refer the reader to
\cite{Sch} for further discussions. A first sufficiently general
solution of this kind of moduli problems was given by the author
in the case that $X$ is a projective curve, $G=\GL(r)$, and
$\rho$ is a homogeneous representation \cite{Sch}. Later, G\'omez
and Sols \cite{GS} established this case on higher dimensional
base manifolds $X$.
\par
The aim of the present paper is to extend these results to the
case when the reductive group is a product of
general linear groups, $G=\GL(r_1)\times\cdots\times \GL(r_t)$, $\rho$ belongs to the class of
homogeneous representations (which comprises all irreducible representations), and $X$ is a base manifold of
arbitrary dimension. However, we will not repeat the detailed
constructions of \cite{Sch}, but rather introduce several
non-trivial ``tricks" which will enable us to adapt the proofs in
that paper to the more general situation studied here.
\par
A nice example of a classification problem which can be
formulated in our context is provided by the work of Casnati and
Ekedahl \cite{CE}. Let $X$ be a projective manifold. Then, any
integral Gorenstein cover $\pi\colon Y\lra X$ of degree $4$ can
be obtained from locally free $\Oh_X$-modules $\E$ and $\F$ of rank
3 and 2, respectively, such that $\det(\E)\cong\det(\F)$, and a
section $s\in H^0(X,\F^\vee\otimes
S^2\E)=\Hom(\F,S^2\E)$. The construction is as follows: If
$\pi\colon \P(\E^\vee):={\rm Proj}(S^*\E)\lra X$ is the projection,
then
\begin{eqnarray*}
\Hom\bigl(\pi^*(\det\E)(-4),\pi^*(\F)(-2)\bigr) & = &
\Hom\bigl(\pi^*(\F^\vee\otimes\det(\E)),\Oh_{\P(\E^\vee)}(2)\bigr)
\\
&=& \Hom\bigl(\F^\vee\otimes\det(\E),S^2\E\bigr)
\\
&=&\Hom(\F,S^2\E).
\end{eqnarray*}
Here,
$\F^\vee\otimes\det(\E)\cong\F^\vee\otimes\det(\F)\cong\F$,
because $\F$ has rank 2. Thus, any section $s\in
H^0(X,\F^\vee\otimes S^2\E)$ yields an exact sequence
$$
\begin{CD}
0 @>>> \pi^*\det(\E)(-4) @>s>> \pi^*(\F)(-2) @>>> \Oh_{\P(\E^\vee)}
@>>> \Oh_Y @>>> 0.
\end{CD}
$$
Hence, the moduli problem for degree 4 covers of $X$ is included
in the moduli problem associated with the group
$\GL(3)\times\GL(2)$ and its representation on $\Hom(\C^2,
S^2\C^3)$. Similarly, degree five covers $\rho\colon Y\lra X$
are determined by locally free sheaves $\E$ and $\F$ of rank $5$
and $4$, respectively, and a homomorphism $\phi\colon \E\lra
\bigwedge^2\F\otimes\det(\E)$ \cite{C5}.
\par
Another interesting moduli problem which we will treat with our
methods comes from the representation theory of finite dimensional
algebras (see \cite{KR} and \cite{Ring} for introductions to this
topic): Let $Q=(V,A,t,h)$ be a quiver with vertex set $V=\{\,
v_1,...,v_t\,\}$ and $\ul{\G} =(\G_a,a\in A)$ a collection of
coherent $\Oh_X$-modules on the projective manifold $X$. This
defines a \it twisted path algebra \rm ${\cal B}={\cal
B}(Q,\ul{\G})$ (see \cite{AG} and \cite{GK}). Modules over ${\cal
B}$ can now be described by \it representations of $Q$\rm, i.e.,
tuples $(\E_v,v\in V; f_a, a\in A)$ composed of $\Oh_X$-modules
$\E_v$, $v\in V$, and twisted homomorphisms $f_a\colon
\G_a\otimes\E_{t(a)}\lra \E_{h(a)}$, $a\in A$. Numerous famous
special cases of this construction have been studied in the
literature, such as the Higgs bundles. Recent research has
focussed on more general aspects of this theory: Gothen and King
\cite{GK} have developed the homological algebra of these
representations and \'Alvarez-C\'onsul and Garc\'\i a-Prada
\cite{AG} formulated a semistability concept for the
representations of $Q$ and proved a Kobayashi-Hitchin
correspondence. The semistability concept depends on additional
parameters $\ul{\sigma}=(\sigma_v\in\Z_{>0},v\in V)$ and
$\ul{\chi}=(\chi_v\in \Q[x], v\in V)$ where the polynomials
$\chi_v$ have degree at most $\dim X-1$, $v\in V$, and an ample
line bundle $\Oh_X(1)$ on $X$. For any coherent sheaf ${\cal A}$
on $X$, the Hilbert polynomial w.r.t.\ $\Oh_X(1)$ is denoted by
$P({\cal A})$. We set
$$
P_{\ul{\sigma},\ul{\chi}}(\E_v,v\in V) :=\sum_{v\in V}
\bigl(\sigma_v P(\E_v)-\chi_v\rk\E_v\bigr)
$$
and
$$
{\rk}_{\ul{\sigma}}(\E_v,v\in V):=\sum_{v\in V}\sigma_v\rk\E_v.
$$
A representation $(\E_v,v\in V; f_a, a\in A)$ is then called \it
(semi)stable\rm, if a) the sheaves $\E_v$, $v\in V$, are torsion free and b)
for any collection of saturated subsheaves $\F_v\subset
\E_v$, $v\in V$, (i.e., $\E_v/\F_v$ is again torsion free, $v\in V$)
not all trivial and not all equal to $\E_v$,
such that $f_a(\G_a\otimes \F_{t(a)})\subset \F_{h(a)}$ for all
arrows $a$, one has
$$
\frac{P_{\ul{\sigma},\ul{\chi}}(\F_v,v\in
V)}{\rk_{\ul{\sigma}}(\F_v,v\in
V)}\q(\preceq)\q\frac{P_{\ul{\sigma},\ul{\chi}}(\E_v,v\in
V)}{\rk_{\ul{\sigma}}(\E_v,v\in V)}.
$$
The notation ``$(\preceq)$" means that ``$\prec$" is used for defining
``stable" and ``$\preceq$" for defining ``semista\-ble", and "$\prec$" and "$\preceq$" refer
to the lexicographic ordering of polynomials. Finally,
$(\E_v,v\in V; f_a, a\in A)$ is called \it polystable\rm, if it is
a direct sum of stable representations $(\E^i_v,v\in V; f^i_a,
a\in A)$, $i=1,...,s$, with
$$
\frac{P_{\ul{\sigma},\ul{\chi}}(\E^i_v,v\in
V)}{\rk_{\ul{\sigma}}(\E^i_v,v\in
V)}=\frac{P_{\ul{\sigma},\ul{\chi}}(\E^j_v,v\in
V)}{\rk_{\ul{\sigma}}(\E^j_v,v\in V)},\q \hbox{for all
$i,j=1,...,s$}.
$$
\par
As one of the interesting and important applications of the main result of this paper, we will prove the following
\begin{Thm*}
Fix Hilbert polynomials $\ul{P}=(P_v, v\in V)$, the sheaves
$\ul{\G}$, as well as the parameters $\ul{\sigma}$, and
$\ul{\chi}$ with $\chi_v=\eta_v\cdot \delta$ for some polynomial
$\delta\in\Q[x]$ and rational numbers $\eta_v$, $ v\in V$.
\par
{\rm i)} There exists a quasi-projective moduli space ${\cal
D}:={\cal D}(Q)^{(\ul{\sigma},\ul{\chi})\rm-ss}_{\ul{P},\ul{\G}}$
for polystable representations $(\E_v,v\in V; f_a,a\in A)$ with
$P(\E_v)=P_v$, $v\in V$. The points corresponding to stable
representations form an open subset ${\cal D}^{\rm s}$.
\par
{\rm ii)} There are a vector space ${\Bbb D}$ and a projective
morphism $H\colon {\cal D}\lra {\Bbb D}$, the {\rm generalized
Hitchin map}.
\end{Thm*}
\begin{Rem*}
\'Alvarez-C\'onsul and Garc\'\i a-Prada define semistability w.r.t.\ parameters 
$\ul{\sigma}=(\sigma_v\in \Q_{>0},\allowbreak v\in V)$ and $\ul{\tau}=(\tau_v\in \Q,v\in V)$.
To be precise, for a representation $(\E_v,v\in V; f_a, a\in A)$, they set
\begin{eqnarray*}
\deg_{\ul{\sigma},\ul{\tau}}(\E_v,v\in V) &:=&
\sum_{v\in V}\bigl(\sigma_v\deg(\E_v)-\tau_v\rk(\E_v)\bigr)
\\
\mu_{\ul{\sigma},\ul{\tau}}(\E_v,v\in V)&:=&
\frac{\deg_{\ul{\sigma},\ul{\tau}}(\E_v,v\in V)}{\rk_{\ul{\sigma}}(\E_v,v\in V)}.
\end{eqnarray*}
\'Alvarez-C\'onsul and Garc\'\i a-Prada say that $(\E_v,v\in V; f_a, a\in A)$ is  \it
(semi)stable\rm, if a) the sheaves $\E_v$, $v\in V$, are torsion free and b)
for any collection of saturated subsheaves $\F_v\subset
\E_v$, $v\in V$, 
such that $f_a(\G_a\otimes \F_{t(a)})\subset \F_{h(a)}$ for all
arrows $a$, one has
$$
\mu_{\ul{\sigma},\ul{\tau}}(\F_v,v\in
V)\q
(\le)
\q \mu_{\ul{\sigma},\ul{\tau}}(\E_v,v\in V).
$$
Since multiplying all the parameters with a common positive factor does not alter the semistability
condition, we may assume that the $\sigma_v$, $v\in V$, be positive integers.
We may choose $\delta:=x^{\dim X-1}$ and $\chi_v:=\tau_v\cdot \delta$, $v\in V$.
Then, $\deg_{\ul{\sigma},\ul{\tau}}(\E_v,v\in V)$ is essentially given as the
coefficient of $x^{\dim X-1}/(\dim X-1)!$ in $P_{\ul{\sigma},\ul{\chi}}$ (see \cite{HL}, Definition 1.2.11, p.\ 13,
for the prescise statement). In particular, we see:
\begin{center}
\begin{tabular}{cl}
&$(\E_v,v\in V; f_a, a\in A)$ stable in the sense \'Alvarez-C\'onsul and Garc\'\i a-Prada
\\
$\Longrightarrow$& stable in our sense
\\
$\Longrightarrow$& semistable in our sense
\\
$\Longrightarrow$& semistable in the sense \'Alvarez-C\'onsul and Garc\'\i a-Prada.
\end{tabular}
\end{center}
Therefore, the above theorem gives a quasi-projective moduli space for the stable objects (in either sense)
together with a Gieseker-type compactification.
\end{Rem*}
In the future, we hope to extend the techniques introduced in this paper to treat the case of other
reductive groups. This provides another motivation for studying the more general and abstract moduli
problems introduced here.
\subsection*{Notation}
$X$ will be a fixed projective manifold over the complex numbers,
and $\Oh_X(1)$ a fixed ample line bundle on $X$. For any coherent
sheaf $\E$, $\deg\E$ is the degree of $\E$ w.r.t.\ to $\Oh_X(1)$,
and $P(\E)$ with $P(\E)(l):=\chi(\E(l))$, for all $l\in{\Bbb N}$,
is the Hilbert polynomial of $\E$ w.r.t.\ $\Oh_X(1)$. In order to
avoid excessive occurrences of the symbol ``$^\vee$", we define
$\P(\cal V)$ as the projective bundle of lines in the fibres of
the vector bundle $\cal V$. For any scheme $S$, $h_S$ denotes its
functor of points $T\lra {\rm Mor}(T,S)$.
\par
In the appendix, we have stated two auxiliary results which will
be used on several occasions.
\subsection*{Acknowledgment}
The author wishes to thank Professor Okonek for reference
\cite{CE}. William Fulton, Bal\'azs Szendr\"oi, and the referee pointed out several typos and inconsistencies.
The author acknowledges support by the DFG through a Heisenberg fellowship and through the priority program
``Globale Methoden in der komplexen Geometrie --- Global Methods in Complex Geometry''. 
\section{Background, formal set-up and statement of the main results}
\subsection{Homogeneous representations}
Let $V=\{\, v_1,...,v_t\,\}$ be a finite index set,
$\ul{r}:=(r_v, v\in V)$ a tuple of positive integers, and define
$$
\GL(V,\ul{r}):=\mathop{\sf X}_{v\in V}\GL(\C^{r_v}).
$$
A (finite dimensional, rational) representation $\rho\colon
\GL(V,\ul{r})\lra \GL(A)$ is said to be \it homogeneous (of degree $\alpha$)\rm, if
there is an integer $\alpha$, such that
$$
\rho(z,...,z)=z^\alpha\cdot {\id}_A,\q \hbox{for all } z\in{\C^*}.
$$
\begin{Ex}
\label{normalf1} 
Every irreducible representation is homogeneous.
\par 
\end{Ex}
For any tuple $\ul{\sigma}=(\sigma_1,...,\sigma_t)$ of positive
integers and $a,b,c\in\Z_{\ge 0}$, we define
$$
W(\ul{\sigma},\ul{r}):= \C^{\sum_{i=1}^t \sigma_ir_i}
$$
and the $\GL(V,\ul{r})$-module
$$
W(\ul{\sigma},\ul{r})_{a,b,c}:=
\bigl(W(\ul{\sigma},\ul{r})^{\otimes a}\bigr)^{\oplus b} \otimes
\bigl(\bigwedge^{\sum_{i=1}^t \sigma_ir_i}
W(\ul{\sigma},\ul{r})\bigr)^{\otimes -c}.
$$
The corresponding representation
$$
\rho_{a,b,c}\colon \GL(V,\ul{r})\lra \GL(W(\ul{\sigma},\ul{r})_{a,b,c})
$$
is homogeneous.
\begin{Prop}
\label{normalf2} 
Let $\ul{\sigma}=(\sigma_1,...,\sigma_t)$ be a tuple of positive integers and $\rho\colon
\GL(V,\ul{r})\lra \GL(A)$ a homogeneous representation.  Then, there are non-negative integers $a$, $b$, and $c$,
such that the module $A$ is a direct summand
of the module
$
W(\ul{\sigma},\ul{r})_{a,b,c}.$ 
\end{Prop}
\begin{proof} 
First, there are integers $s_i,t_i$, $i=1,...,k$, such that $A$ is a direct summand of the $\GL(V,\ul{r})$-module
$$
\bigoplus_{i=1}^k W(\ul{\sigma},\ul{r})^{\otimes s_i}\otimes \Hom\bigl(W(\ul{\sigma},\ul{r}),\C\bigl)^{\otimes t_i}.
$$
This is a consequence of \cite{Deligne}, Proposition 3.1 (a), p.\ 40. Since $\rho$ is assumed to be homogeneous,
we have
$$
s_i-t_i= s_j-t_j,\q \hbox{for all } 1\le i\le j\le k.
$$
The assertion follows now from Corollary~1.2 in \cite{Sch}. 
\end{proof}
\begin{Rem}
The tuple $\ul{\sigma}$ will be
a natural parameter in our theory.
\end{Rem}
\begin{Ex}
\label{QuivExt}
Let $Q=(V,A,t,h)$ be a quiver with vertex set $V=\{\,
v_1,...,v_t\,\}$, arrow set $A=\{\, a_1,...,a_n\,\}$, tail map
$t\colon A\lra V$, and head map $h\colon A\lra V$. Fix a dimension
vector $\ul{r}=(r_v,v\in V)$ as well as another tuple
$\ul{\alpha}=(\alpha_v, v\in V)$ of positive integers with
$\alpha_{h(a)}-\alpha_{t(a)}=\alpha_{h(a^\p)}-\alpha_{t(a^\p)}=:\alpha$ for all $a,a^\p\in A$. Then, the
$\GL(V,\ul{r})$-module
$$
\bigoplus_{a\in A} \Hom \bigl((\C^{r_{t(a)}})^{\otimes
\alpha_{t(a)}},(\C^{r_{h(a)}})^{\otimes \alpha_{h(a)}}\bigr)
$$
is homogeneous of degree $\alpha$.
\end{Ex}
\subsection{$V$-split vector spaces}
First, let $W$ be a finite dimensional $\C$-vector space. A \it
weighted flag in $W$ \rm is a pair $(W^\bullet,\ul{\gamma})$ with
$$
W^\bullet\colon\qquad 0\subsetneq W_1\subsetneq\cdots\subsetneq
W_s\subsetneq W
$$
a --- not necessarily complete --- flag in $W$ and
$\ul{\gamma}=(\gamma_1,...,\gamma_{s+1})$ a vector of integers
with $\gamma_1<\cdots <\gamma_{s+1}$.
\begin{Rem}
\label{WF1} In our context, weighted flags arise in the following
way: Let $\la\colon\C^*\lra \GL(W)$ be a one parameter subgroup
and $\chi_1,...,\chi_{s+1}$, $s\ge 0$, the characters of $\C^*$
with non-trivial eigenspace in $W$. Then,
$\chi_i(z)=z^{\gamma_i}$ with $\gamma_i\in\Z$. Let
$W^{\chi_i}\subset W$ be the corresponding eigenspace. We number
the characters in such a way that $\gamma_1<\cdots<\gamma_{s+1}$.
This yields the weight vector $\ul{\gamma}$. The flag $W^\bullet$
is obtained by setting
$$
W_i:=\bigoplus_{j=1}^i W^{\chi_j},\q i=1,...,s.
$$
\end{Rem}
Let $V=\{\,v_1,...,v_t\,\}$ be an index set. A \it $V$-split
vector space \rm is a collection $(W_v,v\in V)$ of vector spaces
indexed by $V$. Note that $V$-split vector spaces form in a
natural way an Abelian category. A \it weighted flag
$(W_\bullet,\ul{\gamma})$ in the $V$-split vector space
$(W_v,v\in V)$ \rm is a pair $(W_\bullet,\ul{\gamma})$ with
$$
W_\bullet\colon 0\subsetneq (W_1^v,v\in V)\subsetneq
\cdots\subsetneq (W_s^v,v\in V) \subsetneq (W_v,v\in W)
$$
a filtration of $(W_v,v\in V)$ by $V$-split subspaces and
$$
\ul{\gamma}= (\gamma_1,...,\gamma_{s+1})
$$
a vector of integers. We have then the equivalent notions of
\begin{itemize}
\item[a.] Tuples $\bigl((\widehat{W}^v_\bullet,\ul{\gamma}^v),v\in V\bigr)$ of weighted
          flags in the $W_v$, $v\in V$,
          $$
          \widehat{W}^v_\bullet\colon\qquad
          0\subsetneq \widehat{W}^v_1\subsetneq\cdots\subsetneq \widehat{W}^v_{s_v}\subsetneq W_v,
          $$
          $$
          \ul{\gamma}^v=(\gamma^v_1,...,\gamma^v_{s_v+1}),\q v\in V.
          $$
          Here, $s_v=0$ is permitted.
\item[b.] Weighted flags $(W_\bullet,\ul{\gamma})$ in the $V$-split vector space
          $(W_v,v\in V)$.
\end{itemize}
Indeed, suppose we are given a tuple as in a. Let
$\gamma_1<\cdots <\gamma_{s+1}$ be the different weights occurring
among the $\gamma^v_{j}$, $v\in V$, $j=1,...,s_{v}+1$. Then, we
define $(W_j^v,v\in V)$ by $W_j^v:= \widehat{W}_{\iota_v(j)}^v$ with
$$
\iota_v(j):=\max\bigl\{\,\iota=1,...,s_v+1\,|\, \gamma^v_\iota\le \gamma_j\,\bigr\}, 
$$
$v\in
V$,  $j=1,...,s+1$. Conversely, given a weighted flag
$(W_\bullet,\ul{\gamma})$ in the $V$-split vector space
$(W_v,v\in V)$, we get a weighted flag
$(\widehat{W}^v_\bullet,\ul{\gamma}^v)$ in $W_v$ by just
projecting onto $W_v$, $v\in V$. These two operations are clearly
inverse to each other.
\subsubsection*{A weight formula}
Next, we fix $\ul{\sigma}=(\sigma_v\in\Z_{>0}, v\in V)$ and set
$M:=\bigoplus_{v\in V} W_v^{\oplus\sigma_v}$. Suppose we are
given a $V$-split vector space $(W_v,v\in V)$ and a weighted flag
$(W_\bullet,\ul{\gamma})$. Define $M_j:=\bigoplus_{v\in V}
(W_j^v)^{\oplus\sigma_v}$ in order to obtain a weighted flag
$(M_\bullet,\ul{\gamma})$ in $M$.
Assume, furthermore, that we are given quotients $k_v\colon
W_v\lra\C^{t_v}$, $v\in V$, and set $k:=\bigoplus_{v\in
V}k_v^{\oplus\sigma_v}\colon M\lra \C^{t}$, $t:=\sum_{v\in
V}\sigma_v\cdot t_v$. The central formula we will need later is
contained in
\begin{Prop}
\label{weighti} Suppose that, in the above situation, we are
given a tuple $\bigl((\widehat{W}^v_\bullet,\ul{\gamma}^v),v\in
V\bigr)$ of weighted flags in the $W_v$, $v\in V$. Let
$(W_\bullet,\ul{\gamma})$ be the corresponding weighted flag in
$(W_v,v\in V)$ and $(M_\bullet,\ul{\gamma})$ the resulting
weighted flag in $M$. Then, the following identity holds true:
\begin{eqnarray*}
&&\sum_{j=1}^s \frac{\gamma_{j+1}-\gamma_j}{r} \Bigl(r\cdot \dim
k(M_j)-t\cdot \dim M_j\Bigr)\q
\\
&=& \sum_{v\in V} \sigma_v
\biggl(\sum_{j=1}^{s_v}\frac{\gamma^v_{j+1}-\gamma^v_j}{r_v}
\Bigl(r_v\cdot \dim k_v(\widehat{W}_j^v)-t_v\cdot \dim
\widehat{W}_j^v\Bigr)\biggr) \\
&& -\sum_{v\in V}
\sigma_v\cdot\left(\frac{t_v}{r_v}-\frac{t}{r}\right)
\cdot\left(\sum_{j=1}^{s_v+1}\gamma^v_j\bigl(\dim
\widehat{W}^v_j-\dim \widehat{W}^v_{j-1}\bigr)\right).
\end{eqnarray*}
Here, $r_v:=\dim W_v$, $v\in V$, and $r=\sum_{v\in
V}\sigma_v\cdot r_v$.
\end{Prop}
\begin{proof}
From the definitions, the formula
$$
\sum_{v\in V}\sigma_v\cdot\Bigl(
\sum_{j=1}^{s_v+1}\gamma^v_j\bigl(\dim \widehat{W}^v_j-\dim
\widehat{W}^v_{j-1}\bigr) \Bigr) =
\sum_{j=1}^{s+1}\gamma_j\bigl(\dim M_j-\dim M_{j-1}\bigr)
$$
follows immediately. Therefore, the assertion is equivalent to
the following equation
\begin{eqnarray*}
&&\sum_{j=1}^s \frac{\gamma_{j+1}-\gamma_j}{r} \Bigl(r\cdot \dim
k(M_j)-t\cdot \dim M_j\Bigr)
\\
&& -\frac{t}{r}\cdot \sum_{j=1}^{s+1} \gamma_j\cdot(\dim M_j-\dim
M_{j-1})
\\
&=& \sum_{v\in V} \sigma_v
\biggl(\sum_{j=1}^{s_v}\frac{\gamma^v_{j+1}-\gamma^v_j}{r_v}
\Bigl(r_v\cdot \dim k_v(\widehat{W}_j^v)-t_v\cdot \dim
\widehat{W}_j^v\Bigr)\biggr) \\
&& -\sum_{v\in V} \sigma_v\cdot\frac{t_v}{r_v}
\cdot\left(\sum_{j=1}^{s_v+1}\gamma^v_j\bigl(\dim
\widehat{W}^v_j-\dim \widehat{W}^v_{j-1}\bigr)\right).
\end{eqnarray*}
Now,
\begin{eqnarray*}
&&-\frac{t}{r}\cdot \sum_{j=1}^s (\gamma_{j+1}-\gamma_j)\dim M_j
\\
&=& \left(\frac{t}{r}\sum_{j=1}^s \gamma_j \dim M_j\right) -
\left(\frac{t}{r}\sum_{j=1}^s \gamma_{j+1} \dim M_j\right)
\\
&\stackrel{M_0=\{0\}}{=}& \left(\frac{t}{r}\sum_{j=1}^s \gamma_j
\dim M_j\right) - \left(\frac{t}{r}\sum_{j=1}^{s+1} \gamma_{j}
\dim M_{j-1}\right)
\\
&=& -\frac{t}{r}\cdot \gamma_{s+1}\dim
M+\sum_{j=1}^{s+1}\gamma_j\bigl(\dim M_j-\dim M_{j-1}\bigr).
\end{eqnarray*}
Therefore, the left hand side simplifies to
$$
-\frac{t}{r}\cdot \gamma_{s+1}\dim
M+\sum_{j=1}^s(\gamma_{j+1}-\gamma_j)\cdot \dim k(M_j).
$$
With the same argument as before, we see
$$
\sum_{j=1}^s(\gamma_{j+1}-\gamma_j)\cdot \dim
k(M_j)=\gamma_{s+1}\dim k(M)-\sum_{j=1}^{s+1}\gamma_j\bigl(\dim
k(M_j)-\dim k(M_{j-1})\bigr).
$$
Since $\dim M =r$ and $\dim k(M)=t$, the left hand side finally
takes the form
\begin{equation}
\label{LHS} -\sum_{j=1}^{s+1}\gamma_j\bigl(\dim k(M_j)-\dim
k(M_{j-1})\bigr).
\end{equation}
 Likewise, the right hand
side becomes
\begin{equation}
\label{RHS} -\sum_{v\in V}\sigma_v\cdot
\Bigl(\sum_{j=1}^{s_v+1}\gamma^v_j\bigl(\dim
k_v(\widehat{W}^v_j)-\dim k(\widehat{W}^v_{j-1})\bigr)\Bigr).
\end{equation}
The equality of (\ref{LHS}) and (\ref{RHS}) is now clear from the
definitions.
\end{proof}
\begin{Rem}
The conceptual way to see the above formula which explains how
it will arise later is the following: Denote by ${\Bbb G}_v$ the
Gra\ss mannian of $t_v$-dimensional quotients of $W_v$, $v\in V$.
Let ${\frak k}_v\colon W_v\otimes\Oh_{\Bbb G_v}\lra Q_v$ be the
universal quotient and $\Oh_{\Bbb G_v}(1):=\det(Q_v)$, $v\in V$.
Likewise, we let ${\Bbb G}$ be the Gra\ss mannian of
$t$-dimensional quotients of $M$, ${\frak k}\colon
M\otimes\Oh_{\Bbb G}\lra Q$ the universal quotient, and $\Oh_{\Bbb
G}(1):=\det(Q)$. On $\mathop{\sf X}_{v\in V}{\Bbb G}_v$, we have
the quotient
$$
{\frak k}^\p:=\bigoplus_{v\in V}\pi_{\Bbb G_v}^*{\frak k}^{\oplus
\sigma_v}_v \colon M\otimes\Oh_{\mathop{\sf X}_{v\in V}{\Bbb
G}_v}\lra \bigoplus_{v\in V} \pi_{\Bbb G_v}^*Q_v^{\oplus\sigma_v}.
$$
This quotient defines a $\bigl(\mathop{\sf X}_{v\in
V}\GL(W_v)\bigr)$-equivariant embedding
$$
h\colon \mathop{\sf X}_{v\in V}{\Bbb G}_v\hookrightarrow {\Bbb
G},
$$
such that
$$
h^*\Oh_{\Bbb G}(1)=\Oh_{\mathop{\sf X}_{v\in V}{\Bbb
G}_v}(\sigma_{v_1},...,\sigma_{v_t}).
$$
Let $\la_v\colon\C^*\lra\GL(W_v)$ be a one parameter subgroup
which induces the weighted flag
$(\widehat{W}^v_\bullet,\ul{\gamma}^v)$ in $W_v$, $v\in V$, see
Remark~\ref{WF1}. Then, the expression (\ref{RHS}) is just
$$
\mu\bigl(\la,(k_{v_1},...,k_{v_t})\bigr),\q \la:=
(\la_{v_1},...,\la_{v_t}),
$$
w.r.t.\ the linearization of the $\GL(V,\ul{r})$-action in
$\Oh_{\mathop{\sf X}_{v\in V}{\Bbb
G}_v}(\sigma_{v_1},...,\sigma_{v_t})$. Now, we can view $\la$ as
a one parameter subgroup of $\GL(M)$. The induced weighted
filtration of $M$ is then $(M_\bullet,\ul{\gamma})$, and
(\ref{LHS}) agrees with
$$
\mu\bigl(\la, k\bigr)
$$
w.r.t.\ the linearization in $\Oh_{\Bbb G}(1)$. Obviously, we must
have
$$
\mu\bigl(\la, k\bigr)=\mu\bigl(\la,(k_{v_1},...,k_{v_t})\bigr).
$$
\end{Rem}
\subsection{$V$-split sheaves}
\label{split}
We fix a finite index set $V=\{\,v_1,...,v_t\,\}$. A \it
$V$-split sheaf \rm is simply a tuple $(\E_v, v\in V)$ of coherent
sheaves on $X$. Likewise, a \it homomorphism between $V$-split
sheaves \rm $(\E_v,v\in V)$ and $(\E^\p_v, v\in V)$ is a
collection $(f_v,v\in V)$ of homomorphisms $f_v\colon  \E_v\lra
\E_v^\p$, $v\in V$. In this way, the $V$-split sheaves on $X$
form an Abelian category. The \it type of the $V$-split sheaf
$(\E_v,v\in V)$ \rm is the tuple $\ul{P}(\E_v,v\in
V)=(P(\E_{v_1}),...,P(\E_{v_t}))$.
\begin{Rem}
The datum of a $V$-split vector bundle is equivalent to the datum 
of a principal $\GL(V,\ul{r})$-bundle. Thus, a $V$-split sheaf can be
seen as the natural "singular" version of a principal $\GL(V,\ul{r})$-bundle.
\end{Rem}
Now, we fix additional parameters
$\ul{\sigma}=(\sigma_v\in\Z_{>0},v\in V)$ and
$\ul{\chi}=(\chi_v\in \Q[x], v\in V)$ where the polynomials
$\chi_v$ have degree at most $\dim X-1$, $v\in V$. We denote by
$\ol{\chi}_v$ the coefficient of $x^{\dim X-1}$ in $\chi_v$,
$v\in V$. Then, we define for any $V$-split sheaf $(\E_v,v\in V)$:
\par
\begin{tabular}{ll}
\\
The \it $(\ul{\sigma},\ul{\chi})$-degree\rm: & $\deg_{\ul{\sigma},\ul{\chi}}(\E_v,v\in V)
:=\sum_{v\in V}
\bigl(\sigma_v\deg\E_v-\ol{\chi}_v\rk\E_v\bigr)$;
\\
\\
The \it $(\ul{\sigma},\ul{\chi})$-Hilbert polynomial\rm: &
$P_{\ul{\sigma},\ul{\chi}}(\E_v,v\in V)
:=\sum_{v\in V}
\bigl(\sigma_v P(\E_v)-\chi_v\rk\E_v\bigr)$;
\\
\\
The \it $\ul{\sigma}$-rank\rm: & $\rk_{\ul{\sigma}}(\E_v,v\in V):=\sum_{v\in V}\sigma_v\rk\E_v$;
\\
\\
The \it $(\ul{\sigma},\ul{\chi})$-slope\rm: & $\mu_{\ul{\sigma},\ul{\chi}}(\E_v,v\in V)
:=\frac{\deg_{\ul{\sigma},\ul{\chi}}(\E_v,v\in V)}{\rk_{\ul{\sigma}}(\E_v,v\in V)}$.
\\
\\
\end{tabular}
\par
Note that the $(\ul{\sigma},\ul{\chi})$-degree, $(\ul{\sigma},\ul{\chi})$-Hilbert-polynomial,
and $\ul{\sigma}$-rank all behave additively on short exact sequences. Thus, the
$(\ul{\sigma},\ul{\chi})$-slope will have all the formal properties of the usual slope.
More specifically, we call a $V$-split sheaf $(\E_v,v\in V)$ \it
$(\ul{\sigma},\ul{\chi})$-slope (semi)stable\rm, if $\E_v$ is torsion free, $v\in V$, and,
for all non-trivial, proper
$V$-split subsheaves $(\F_v,v\in V)$, the inequality
$$
\mu_{\ul{\sigma},\ul{\chi}}\bigl(\F_v,v\in V\bigr)(\le)
\mu_{\ul{\sigma},\ul{\chi}}\bigl(\E_v,v\in V\bigr)
$$
is satisfied.
\begin{Rem}
Let $(\E_v,v\in V)$ be a $(\ul{\sigma},\ul{\chi})$-slope semistable $V$-split sheaf.
Let $v_0$ be a vertex with $\E_{v_0}\neq 0$
and look at a subsheaf $0\subsetneq \F\subseteq \E_{v_0}$.
Set $\F_v:= 0$ for $v\neq v_0$ and $\F_{v_0}:=\F$. Then, we get
$\mu_{\ul{\sigma},\ul{\chi}}\bigl(\F_v,v\in V\bigr)=\mu(\F)-(\ol{\chi}_{v_0}/\sigma_{v_0})$,
so that the semistability condition yields
$$
\mu(\F)\le \mu_{\ul{\sigma},\ul{\chi}}\bigl(\E_v,v\in V\bigr)+
\frac{\ol{\chi}_{v_0}}{\sigma_{v_0}}.
$$
Likewise, we find for every non-trivial quotient $\cal Q$ of $\E_{v_0}$
$$
\mu(\cal Q)
\ge
\mu_{\ul{\sigma},\ul{\chi}}\bigl(\E_v,v\in V\bigr)+
\frac{\ol{\chi}_{v_0}}{\sigma_{v_0}}.
$$
If we apply this to $\F=\E_{v_0}=\cal Q$, we see that $\E_{v_0}$ must be
a semistable sheaf with slope $
\mu_{\ul{\sigma},\ul{\chi}}\bigl(\E_v,v\in V\bigr)+
{\ol{\chi}_{v_0}}/{\sigma_{v_0}}$. Note that this forces $\sum_{v\in V} r_v\cdot \ol{\chi}_v
=0$.
\end{Rem}
Recall that any homomorphism $f\colon \E\lra\E^\p$ between the semistable sheaves
$\E$ and $\E^\p$ will be zero, if $\mu(\E)>\mu(\E^\p)$.
Therefore, we deduce
\begin{Prop}
\label{homtriv}
Let $(\E_v,v\in V)$ and $(\E^\p_v,v\in V)$ be
$(\ul{\sigma},\ul{\chi})$-slope semistable $V$-split sheaves. Assume
$$
\mu_{\ul{\sigma},\ul{\chi}}\bigl(\E_v,v\in V\bigr)
+\min_{v\in V}\left\{\,\frac{\ol{\chi}_{v}}{\sigma_{v}}\,\right\}
>
\mu_{\ul{\sigma},\ul{\chi}}\bigl(\E^\p_v,v\in V\bigr)
+\max_{v\in V}\left\{\,\frac{\ol{\chi}_{v}}{\sigma_{v}}\,\right\}.
$$
Then, for any choice of exponents $s_v, s^\p_v>0$, $v\in V$,
any homomorphism $f\colon  \bigoplus_{v\in V} \E_v^{\oplus s_v}\lra
\bigoplus_{v\in V}{\E^\p_v}^{\oplus s^\p_v}$
is zero.
\end{Prop}
Finally, we have
\begin{Prop}[Harder-Narasimhan filtration]
\label{HS}
Let $(\E_v,v\in V)$ be any torsion free $V$-split sheaf. Then, $(\E_v,v\in V)$
possesses a unique filtration by $(\ul{\sigma},\ul{\chi})$-destabilizing
$V$-split subsheaves
$$
0=(\F_0^v,v\in V)\subsetneq (\F_1^v,v\in V)\subsetneq\cdots\subsetneq (\F_s^v,v\in V)
\subsetneq (\F_{s+1}^v,v\in V)=(\E_v,v\in V),
$$
such that
\begin{enumerate}
\item[\rm 1.] The $V$-split sheaf
$(\F_i^v/\F_{i-1}^v,v\in V)$ is $(\ul{\sigma},\ul{\chi})$-slope semistable,
$i=1,...,s+1$.
\item[\rm 2.]
$
\mu_{\ul{\sigma},\ul{\chi}}\bigl(\F_i^v/\F_{i-1}^v,v\in V\bigr)
>
\mu_{\ul{\sigma},\ul{\chi}}\bigl(\F_{i+1}^v/\F_{i}^v,v\in
V\bigr)$, $i=1,...,s$.
\end{enumerate}
\end{Prop}
\subsubsection*{The weight formula for sheaves}
Let $(\E_v,v\in V)$ be a $V$-split sheaf. As before, the
following data are equivalent
\begin{itemize}
\item[a.] Tuples $\bigl((\widehat{\E}^v_\bullet,\ul{\gamma}^v),v\in V\bigr)$ of weighted
          filtrations of the $\E_v$, $v\in V$,
          $$
          \widehat{\E}^v_\bullet\colon\qquad
          0\subsetneq \widehat{\E}^v_1\subsetneq\cdots\subsetneq \widehat{\E}^v_{s_v}\subsetneq \E_v,
          $$
          $$
          \ul{\gamma}^v= (\gamma^v_1,...,\gamma^v_{s_v+1}),\q v\in V.
          $$
\item[b.] Weighted filtrations $(\E_\bullet,\ul{\gamma})$ of the $V$-split sheaf
          $(\E_v,v\in V)$.
\end{itemize}
Moreover, given $\ul{\sigma}=(\sigma_v\in\Z_{>0},v\in V)$, a
$V$-split sheaf $(\E_v,v\in V)$, and a weighted filtration
$(\E_\bullet,\ul{\gamma})$, we set $\E^{\rm
total}:=\bigoplus_{v\in V}\E_v^{\oplus\sigma_v}$ and
$$
\E^{\rm total}_\bullet\colon\qquad 0\subsetneq \E_1^{\rm total}:=
\bigoplus_{v\in
V}{\E_1}^{v,\oplus\sigma_v}\subsetneq\cdots\subsetneq \E_s^{\rm
total}:=\bigoplus_{v\in V}{\E_s}^{v,\oplus\sigma_v} \subsetneq
\E^{\rm total}.
$$
\begin{Prop}
\label{weightii} Suppose that, in the above situation, we are
given a tuple $\bigl((\widehat{\E}^v_\bullet,\ul{\gamma}^v),v\in
V\bigr)$ of weighted filtrations of the $\E_v$, $v\in V$. Let
$(\E^{\rm total}_\bullet,\ul{\gamma})$ be the resulting weighted
filtration of $\E^{\rm total}$. Then, for all $l\gg 0$,
\begin{eqnarray*}
&&\sum_{j=1}^s \frac{\gamma_{i+1}-\gamma_i}{P(\E^{\rm total})(l)}
\Bigl(P(\E^{\rm total})(l)\cdot \rk\E_j^{\rm total}-P(\E_j^{\rm
total})(l)\cdot \rk\E^{\rm total} \Bigr)\
\\
&=&
\sum_{v\in V} \sigma_v
\biggl(\sum_{j=1}^{s_v}\frac{\gamma^v_{j+1}-\gamma^v_j}{P(\E_v)(l)}
\Bigl(P(\E_v)(l)\cdot \rk \widehat{\E}^v_j-
P(\widehat{\E}^v_j)(l)\cdot \rk\E_v\Bigr)\biggr)
\\
&& -\sum_{v\in V}
\sigma_v\cdot\left(\frac{\rk\E_v}{P(\E_v)(l)}-\frac{\rk\E^{\rm
total}}{P(\E^{\rm total})(l)}\right)
\cdot\left(\sum_{j=1}^{s_v+1}\gamma^v_j\bigl(
P(\widehat{\E}^v_j)(l)-P(\widehat{\E}^v_{j-1})(l)\bigr)\right).
\end{eqnarray*}
\end{Prop}
\noindent\it Proof\rm. For $l\gg 0$, we have
\begin{itemize}
\item $\widehat{\E}^v_j(l)$ is globally generated, $v\in V$, $j=1,...,s_v+1$,
\item $H^i(\widehat{\E}^v_j(l))=0$, $i>0$, $v\in V$, and $j=1,...,s_v+1$.
\end{itemize}
Then, we may write $\E_v(l)$ as a quotient $q_v\colon
\Oh_X^{\oplus P(\E_v)(l)} \lra \E_v(l)$, such that $H^0(q_v)$ is
an isomorphism, $v\in V$. Restricting this to a general point
$x\in X$ yields $k_v\colon \C^{P(\E_v)(l)}\lra \C^{\rk\E_v}$,
$v\in V$. Now, apply Proposition~\ref{weighti} to the tuple
$\Bigl(\bigl(H^0(\widehat{\E}^v_\bullet(l)),\ul{\gamma}^v\bigr),
v\in V\Bigr)$ under the identification of $\C^{P(\E_v)(l)}$ with
$H^0(\E_v(l))$. \qed
\subsection{$\rho$-pairs}
In this section, we will fix a dimension vector $\ul{r}=(r_v,
v\in V)$ and a homogeneous representation $\rho\colon \GL(V,\ul{r})\lra
\GL(A)$. In slight deviation from the conventions in the
introduction (see Remark~\ref{explain}), the objects we would like
to consider are pairs $(\E_v, v\in V;\Psi)$ where $(\E_v, v\in V)$
is a $V$-split vector bundle, such that $\rk \E_v=r_v$, $v\in V$,
and $\Psi\colon X\lra \P(\F^\vee_\rho)$ is a section. Here,
$\F_\rho$ is the vector bundle with fibre $A$ associated to
$(\E_v, v\in V)$ via the representation $\rho$. Now, the section
$\Psi$ is specified by a line bundle $\N$ and a surjective
homomorphism $\phi\colon \F_\rho\lra \N$, and two such
homomorphisms will yield the same section, if and only if they
differ by a constant $z\in\C^*$. Thus, in order to find
projective moduli spaces (at least over curves), we consider
tuples $(\E_v, v\in V;\N;\phi)$ where $\N$ is a line bundle and
$\phi\colon \F_\rho\lra \N$ is a non-trivial homomorphism. Such an
object will be referred to as a \it $\rho$-pair\rm, and the tuple
$(\ul{P}(\E_v, v\in V),\N)$ is called the \it type\rm. Two
$\rho$-pairs $(\E_v, v\in V;\N;\phi)$ and $(\E^\p_v, v\in
V;\N^\p;\phi^\p)$ will be considered \it equivalent\rm, if there
are isomorphisms $\psi_v\colon \E_v\lra \E^\p_v$, $v\in V$, and
$z\colon \N\lra \N^\p$, such that
$$
\phi^\p=z\circ\phi\circ \psi_\rho^{-1},
$$
with $\psi_\rho \colon \F_\rho\lra \F^\p_\rho$ the isomorphism
induced by the $\psi_v$, $v\in V$.
\par
Given a tuple $\ul{P}=(P_v, v\in V)$ of Hilbert polynomials and a
line bundle $\n$, a \it family of $\rho$-pairs of type
$(\ul{P},\n)$ parameterized by the scheme $S$ \rm is the datum of a
tuple $(\E_{S,v}, v\in V; \EL_S, \phi_S)$ with vector bundles
$\E_{S,v}$, $v\in V$, on $S\times X$, such that
$P(\E_{S,v|\{s\}\times X})=P_v$ for all $s\in S$, $v\in V$,
$\EL_S$ a line bundle on $S$, and $\phi_S\colon \F_{S,\rho}\lra
\pi_S^*\EL_S\otimes \pi_X^*\n$ a homomorphism with non trivial
restriction to every fibre $\{s\}\times X$, $s\in S$. Two such
families $(\E_{S,v}, v\in V; \EL_S, \phi_S)$ and $(\E^\p_{S,v},
v\in V; \EL^\p_S, \phi^\p_S)$ will be considered \it
equivalent\rm, if there are isomorphisms $\psi_{S,v}\colon
\E_{S,v}\lra \E_{S,v}^\p$, $v\in V$, and $z_S\colon \EL_S\lra
\EL_S^\p$, such that
$$
\phi_S^\p= \bigl(\pi_S^*(z_S)\otimes
{\id}_{\pi_X^*\n}\bigr)\circ \phi_S\circ\psi_{S,\rho}^{-1}.
$$
\begin{Rem}
\label{explain} First, we note that using $\F_\rho^\vee$ instead
of $\F_\rho$ is for notational convenience only. Then, the right
analogue to the problems mentioned in the introduction would be
the study of tuples $(\E_v, v\in V; \phi)$ where $\phi\colon
\F_\rho\lra\Oh_X$ is a surjective homomorphism, and $(\E_v, v\in V;
\phi)$ and $(\E^\p_v, v\in V; \phi^\p)$ should be identified, if
and only if there are isomorphisms $\psi_v\colon \E_v\lra \E_v^\p$
with $\phi^\p=\phi\circ \psi_\rho^{-1}$. In the case of a homogeneous
representation of non-zero degree, this
equivalence relation will identify $(\E_v, v\in V; \phi)$ and
$(\E_v, v\in V; z\cdot \phi)$, $z\in\C^*$, anyway.
Otherwise, one may add the trivial representation to $\rho$. This
means that we consider tuples $(\E_v, v\in V; \phi, \eps)$ with
$(\E_v, v\in V; \phi)$ as before and $\eps\in\C$, but the
equivalence relation becomes $(\E_v, v\in V; \phi, \eps)\sim
(\E^\p_v, v\in V; \phi^\p, \eps^\p)$, if there are isomorphisms
$\psi_v\colon \E_v\lra \E_v^\p$, $v\in V$, and a $z\in\C^*$, such
that
$$
z\cdot \phi^\p =\phi\circ \psi_\rho^{-1},\q \hbox{and}\
z\cdot\eps^\p=\eps.
$$
Then, we may recover the original objects in the form $(\E_v, v\in
V; \phi, 1)$. Thus, our concept is more flexible rather than more
restrictive than the one presented in the introduction.
\end{Rem}
In order to define the semistability concept we introduce
additional parameters:
\begin{itemize}
\item a tuple $\ul{\sigma}=(\sigma_v, v\in V)$ of positive
integers.
\item a tuple $\ul{\eta}=(\eta_v, v\in V)$ of rational numbers,
such that $\sum_{v\in V}\eta_vr_v=0$.
\item a positive polynomial $\delta\in \Q[x]$. Then,
$\chi_v:=\eta_v\cdot \delta$, $v\in V$.
\end{itemize}
\par
Given any torsion free $\Oh_X$-module $\E$, we call a submodule
$\F\subset\E$ \it saturated\rm, if the quotient $\E/\F$ is still
torsion free. The test objects for the semistability concept will
be weighted filtrations $(\E_\bullet,\ul{\gamma})$ of the
$V$-split vector bundle $(\E_v,v\in V)$ where each $(\E_j^v,v\in
V)$ consists of saturated subsheaves $\E_j^v\subset\E_v$, $v\in
V$, $j=1,...,s$. For such a weighted filtration, we define
$\ul{\alpha}=(\alpha_1,...,\alpha_{s+1})$ by
$\alpha_i:=(\gamma_{i+1}-\gamma_i)/R$, $R:=\rk(\bigoplus_{v\in V}
\E_v^{\oplus\sigma_v})$, $i=1,...,s$. We now set
$$
M_{\ul{\sigma},\ul{\chi}}(\E_\bullet,\ul{\alpha})\ :=\
\sum_{j=1}^s\alpha_j\cdot
\Bigl(P_{\ul{\sigma},\ul{\chi}}(\E_v,v\in
V)\cdot{\rk}_{\ul{\sigma}}(\E_j^v,v\in V) -
P_{\ul{\sigma},\ul{\chi}}(\E_j^v,v\in V)\cdot
{\rk}_{\ul{\sigma}}(\E_v,v\in V)\Bigr).
$$
If we are also given a non-trivial homomorphism $\phi\colon
\F_\rho\lra\N$, we have to define the quantity
$$
\mu(\E_\bullet,\ul{\gamma};\phi).
$$
The natural, though complicated, definition for $\#V=1$ was
explained in \cite{Sch}. We adapt it to our setting. Let
$M:=\bigoplus_{v\in V} \C^{\sigma_v\cdot r_v}$. Then, for
appropriate $a,b,c$, the module $A$ will be a submodule of
${M^{\otimes a}}^{\oplus b}\otimes (\bigwedge^{\dim W}
W)^{\otimes -c}$. We introduce $\E^{\rm total}_j:=\bigoplus_{v\in
V}(\E_j^v)^{\oplus\sigma_v}$, $j=1,...,s$. There is the weighted
filtration $(\E_\bullet^{\rm total},\ul{\alpha})$
$$
0\subsetneq \E_1^{\rm total}\subsetneq\cdots\subsetneq\E_s^{\rm
total}\subsetneq \E^{\rm total}
$$
of $\E^{\rm total}=\bigoplus_{v\in V}\E_v^{\oplus\sigma_v}$. We
choose a flag
$$
0\subsetneq M_1\subsetneq\cdots\subsetneq M_s\subsetneq M
$$
with $\dim M_i=\rk \E^{\rm total}_i$, $i=1,...,s$. Over a
suitable open subset ${\cal U}$, the homomorphism $\phi$ will be
surjective, and there will be a trivialization $\psi\colon
M\otimes\Oh_{\cal U}\cong \E^{\rm total}$ with
$\psi(M_i\otimes\Oh_{\cal U})=\E^{\rm total}_i$, $i=1,...,s$. We
may write
$$
\F_{a,b,c}:={{\E^{\rm total}}^{\otimes a}}^{\oplus b}\otimes
\det\bigl(\E^{\rm total}\bigr)^{\otimes -c}=\F_\rho\oplus
\F^\p,
$$
so that $\phi$ and the trivialization $\psi$ yield a morphism
$$
\sigma\colon {\cal U}\lra \P(\F^\vee_\rho)\hookrightarrow
\P(\F_{a,b,c}^\vee) \cong \P\bigl( {{M^{\otimes a}}^{\oplus
b}}^\vee\bigr)\times {\cal U}\lra \P\bigl( {{M^{\otimes
a}}^{\oplus b}}^\vee\bigr).
$$
After the choice of a one parameter subgroup $\la\colon \C^*\lra
\SL(M)$ which induces the weighted filtration
$(M_\bullet,\ul{\gamma}^\p)$ with
$$
\ul{\gamma}^\p=(\gamma_1^\p,...,\gamma_{s+1}^\p)=\sum_{i=1}^s\alpha_i
(\underbrace{R_i-R,...,R_i-R}_{R_i\times},\underbrace{R_i,...,R_i}_{(R-R_i)\times}),\q
R_i:=\rk \E_i^{\rm total}, i=1,...,s,
$$
we set
$$
\mu(\E_\bullet,\ul{\gamma};\phi):={\rm max}\bigl\{\,
\mu(\la,\sigma(x))\,|\, x\in {\cal U}\,\bigr\}.
$$
As in \cite{Sch}, one verifies that this is, in fact, well
defined.
\begin{Rem}
\label{explain3} i) One might define
$\mu(\E_\bullet,\ul{\gamma};\phi)$ without reference to the
embedding of $A$ into ${M^{\otimes a}}^{\oplus b}\otimes
(\bigwedge^{\dim W} W)^{\otimes -c}$, by working with one
parameter subgroups of $\GL(V,\ul{r})$.
\par
ii) An easier, more elegant, and equivalent definition \cite{GS} is
$$
\mu(\E_\bullet,\ul{\gamma};\phi):=-\min_{(i_1,...,i_a)\in\atop
\{\,1,...,s+1\,\}^{\times a}}\Bigl\{\,
\gamma^\p_{i_1}+\cdots+\gamma^\p_{i_a}\,|\, \hbox{$\tau$ rest.\
to $(\E^{\rm total}_{i_1}\otimes\cdots\otimes \E^{\rm
total}_{i_a})^{\oplus b}$ is non trivial}\,\Bigr\}.
$$
Here, $\tau\colon\F_{a,b,c}\lra \F_\rho\stackrel{\phi}{\lra}\N$.
However, for the computations in examples, the above definition
turns out to be more useful (see \cite{Sch}).
\end{Rem}
\begin{Conv}
Since the quantities introduced above depend only on
$\ul{\alpha}$, we will refer to a pair $(\E_\bullet,
\ul{\alpha})$, composed of a filtration
$$
\E_\bullet\colon\qquad 0\subsetneq (\E_1^v,v\in V)\subsetneq
\cdots\subsetneq (\E_s^v,v\in V) \subsetneq (\E_v,v\in V)
$$
of $(\E_v,v\in V)$ by non-trivial, proper, and saturated
$V$-split subsheaves and a tuple
$\ul{\alpha}=(\alpha_1,...,\alpha_s)$ of positive rational
numbers, as a \it weighted filtration \rm in the future.
\end{Conv}
A $\rho$-pair $(\E_v,v\in V;\N,\phi)$ will be called \it
$(\ul{\sigma},\ul{\eta},\delta)$-(semi)stable \rm or just \it
(semi)stable\rm, if for every weighted filtration $(\E_\bullet,
\ul{\alpha})$ of $(\E_v,v\in V)$
$$
M_{\ul{\sigma},\ul{\chi}}(\E_\bullet,\ul{\alpha})
+\delta\cdot\mu(\E_\bullet,\ul{\alpha};\phi) (\succeq)  0.
$$
A few comments are in order.
\begin{Rem}
\label{explain2} Since $\sum_{v\in V}\chi_v\cdot r_v\equiv 0$, we
may write
\begin{eqnarray*}
&&M_{\ul{\sigma},\ul{\chi}}(\E_\bullet,\ul{\alpha}) =
\rk\Bigl(\bigoplus_{v\in
V}\E_v^{\oplus\sigma_v}\Bigr)\cdot\sum_{j=1}^s\alpha_j\Bigl(
\sum_{v\in V}\chi_v\rk\E_j^v\Bigr)
\\
&&\q +\q \underbrace{\sum_{j=1}^s\alpha_j\biggl(
P\Bigl(\bigoplus_{v\in V}\E_v^{\oplus\sigma_v}\Bigr)\cdot
\rk\Bigl(\bigoplus_{v\in V}\E_j^{v,\oplus\sigma_v}\Bigr) -
P\Bigl(\bigoplus_{v\in V}\E_j^{v,\oplus\sigma_v}\Bigr)\cdot
\rk\Bigl(\bigoplus_{v\in V}\E_v^{\oplus\sigma_v}\Bigr)
\biggr)}_{=:M_{\ul{\sigma}}(\E_\bullet,\ul{\alpha})}.
\end{eqnarray*}
A bounded family of $V$-split vector bundles $(\E_v, v\in V)$ may
be parameterized by a product of quot schemes $\ul{\frak Q}={\sf
X}_{v\in V} {\frak Q}_v$. Assigning to a point $([q_v], v\in
V)\in \ul{\frak Q}$ the quotient $\bigl[\bigoplus_{v\in V}
q_v^{\oplus \sigma_v}\bigr]$ induces an injective and proper
morphism from $\ul{\frak Q}$ to some other quot scheme
$\widetilde{\frak Q}$. In this way, we can induce linearizations
on $\ul{\frak Q}$ by linearizations on $\widetilde{\frak Q}$, and
this shows how the quantity
$M_{\ul{\sigma}}(\E_\bullet,\ul{\alpha})$ is obtained.
The linearization of the $\GL(V,\ul{r})$-action on the space
$\P({M^{\otimes a}}^{\oplus b})$ in $\Oh(1)$ induced by
$\rho_{a,b,c}$ may be modified by a character, so that the
determinant on ${M^{\otimes a}}^{\oplus b}$ induces the trivial
character on the center ${\cal Z}$ of $\GL(V,\ul{r})$, and the
quantity $\mu(\E_\bullet,\ul{\alpha};\phi)$ has been defined
w.r.t.\ such a linearization. The parameter $\delta$ reflects the
fact that the given linearization in $\Oh(1)$ may be raised to
some tensor power.
Finally, any linearization might be altered by a character $\chi$
of $\GL(V,\ul{r})$. The choice of such a character
is encoded by the rational numbers $\eta_v$, $v\in V$.
These considerations explain how the semistability concept we have
introduced naturally  ``mixes" the semistability concept for
vector bundles and the invariant theory of the representation
$\rho$.
The condition $\sum_{v\in V}\eta_v r_v=0$ is used to simplify some computations. 
It can, however, be assumed \sl without loss of generality\rm.
For this, note that for arbitrary parameters $\sigma_v$ and $\eta_v$, $v\in V$, 
\begin{itemize}
\item $M_{\ul{\sigma},\ul{\chi}}(\E_\bullet,\ul{\alpha})$ is defined the same way,
\item $\mu(\E_\bullet,\ul{\gamma};\phi)$ does not depend on the $\eta_v$, $v\in V$.
\end{itemize}
In particular, we may define (semi)stability w.r.t.\ these parameters.
Suppose we are given arbitrary rational numbers $\eta_v$, $v\in V$, and $d\in \Q$.
Define $\eta^\p_v:=\eta_v-d\cdot \sigma_v$, $v\in V$. Then, 
$$
P_{\ul{\sigma},\ul{\chi}^\p}(\F_v,v\in V)=P_{\ul{\sigma},\ul{\chi}}(\F_v,v\in V)+
d\cdot\delta\cdot \left(\sum_{v\in V}\sigma_v \rk(\F_v)\right).
$$
It follows that
\begin{eqnarray*}
&& P_{\ul{\sigma},\ul{\chi}}(\E_v,v\in V)\cdot {\rk}_{\ul{\sigma}}(\F_v,v\in V)-
   P_{\ul{\sigma},\ul{\chi}}(\F_v,v\in V)\cdot {\rk}_{\ul{\sigma}}(\E_v,v\in V)
\\
&=& P_{\ul{\sigma},\ul{\chi}^\p}(\E_v,v\in V)\cdot {\rk}_{\ul{\sigma}}(\F_v,v\in V)-
   P_{\ul{\sigma},\ul{\chi}^\p}(\F_v,v\in V)\cdot {\rk}_{\ul{\sigma}}(\E_v,v\in V),
\end{eqnarray*}
so that always
$$
M_{\ul{\sigma},\ul{\chi}}(\E_\bullet,\ul{\alpha})=M_{\ul{\sigma},\ul{\chi}^\p}(\E_\bullet,\ul{\alpha}),
$$
and the concept of (semi)stability defined w.r.t.\ to the parameters $\sigma_v$, $\eta_v$, $v\in V$, and
$\delta$ equals the one defined w.r.t.\ to the parameters $\sigma_v$, $\eta^\p_v$, $v\in V$, and
$\delta$. If we now set
$$
d:=\frac{\sum_{v\in V}\eta_v r_v}{\sum_{v\in V} \sigma_vr_v}, 
$$
then
$$
\sum_{v\in V}\eta_v^\p r_v=\left(\sum_{v\in V}\eta_vr_v\right)-d\cdot \left(\sum_{v\in V} \sigma_vr_v\right)=0.
$$
\end{Rem}
We fix the Hilbert polynomials $\ul{P}=(P_v,v\in V)$ and the line
bundle $\n$ and define the functors
\begin{eqnarray*}
\ul{\rm M}(\rho)^{(\ul{\sigma},\ul{\eta},\delta)-\rm
(s)s}_{\ul{P}/\n}\colon \q \ul{\rm Sch}_\C &\lra &\ul{\rm Sets}
\\
S &\lma &\left\{\,\hbox{Equivalence classes of families of
(semi)stable} \atop \hbox{$\rho$-pairs of type $(\ul{P},\n)$
parameterized by $S$}\,\right\}.
\end{eqnarray*}
\begin{Thm}
\label{main0} {\rm i)} If $\dim(X)=1$, there exist a projective
scheme ${\cal M}:={\cal
M}(\rho)^{(\ul{\sigma},\ul{\eta},\delta)-\rm ss}_{\ul{P}/\n}$ and
a natural transformation $\theta\colon  \ul{\rm
M}(\rho)^{(\ul{\sigma},\ul{\eta},\delta)-\rm ss}_{\ul{P}/\n} \lra
h_{\cal M}$, such that for any other scheme ${\cal M}^\p$ and any
other natural transformation $\theta^\p\colon \ul{\rm
M}(\rho)^{(\ul{\sigma},\ul{\eta},\delta)-\rm ss}_{\ul{P}/\n}\lra
h_{\cal M^\p}$, there is a unique morphism $\zeta\colon \cal
M\lra \cal M^\p$ with $\theta^\p=h_\zeta\circ\theta$. The space
${\cal M}$ contains an open subscheme ${\cal M}^{\rm s}$ which is
a coarse moduli scheme for the functor $\ul{\rm
M}(\rho)^{(\ul{\sigma},\ul{\eta},\delta)-\rm s}_{\ul{P}/\n}$.
\par
{\rm ii)} If $\dim(X)\ge 1$, there exists a quasi-projective
coarse moduli scheme ${\cal
M}(\rho)^{(\ul{\sigma},\ul{\eta},\delta)-\rm s}_{\ul{P}/\n}$ for
the functor $\ul{\rm M}(\rho)^{(\ul{\sigma},\ul{\eta},\delta)-\rm
s}_{\ul{P}/\n}$.
\end{Thm}
\begin{Ex}
The above results yield a semistability concept for coverings in
the description of Casnati and Ekedahl (as reviewed in the
introduction) and provides moduli spaces.
\end{Ex}
It is obvious that the theorem has to be proved only for
representations of the type $\rho_{a,b,c}$. If $\dim(X)>1$, then
in order to compactify the moduli spaces, one needs also torsion
free sheaves. It is however not clear how to associate a sheaf
$\F_\rho$ to a $V$-split torsion free sheaf $(\E_v,v\in V)$ via
an arbitrary representation $\rho$. However, for representations
of the form $\rho_{a,b,c}$, this is obvious and one obtains
natural compactifications. In the setting of quiver
representations, we will exhibit another natural method to reduce
a moduli problem to one for $\rho_{a,b,c}$-pairs. This
illustrates the importance and the usefulness of the theory which
we will outline in the next section.
\subsection{Decorated $V$-split sheaves}
\label{deko1}
We fix the following data
\begin{itemize}
\item a tuple of Hilbert polynomials $\ul{P}=(P_v,v\in V)$;
\item a positive polynomial $\delta\in\Q[x]$ of degree at most $\dim X-1$;
\item a tuple of rational numbers $\ul{\eta}=(\eta_v, v\in V)$ with
$\sum_{v\in V} \eta_v\cdot r_v=0$.
      Here, $r_v$ is the rank dictated by the Hilbert polynomial $P_v$, $v\in V$.
      Define $\chi_v:=\eta_v\cdot\delta$, $v\in V$.
\item a tuple $\ul{\sigma}=(\sigma_v,v\in V)$ of positive integers.
\end{itemize}
Given a $V$-split sheaf $(\E_v,v\in V)$ of type $\ul{P}$ and
non-negative integers $a$, $b$, $c$, and $m$, a \it decoration of
type $(a,b,c,m)$ on $(\E_v,v\in V)$ \rm is a homomorphism
$$
\tau\colon  \biggl(\Bigl(\bigoplus_{v\in V}\E_{v}^{\oplus \sigma_v}\Bigr)^{\otimes a}\biggr)^
{\oplus b}
\lra \det\Bigl(\bigoplus_{v\in V}\E_{v}^{\oplus \sigma_v}\Bigr)^{\otimes c}\otimes\Oh_X(m).
$$
\par
Two tuples $(\E_v,v\in V;\tau)$ and $(\E^\p_v,v\in V;\tau^\p)$
are called \it equivalent\rm, if there are $z\in\C^*$ and isomorphisms $\psi_v
\colon  \E_v\lra \E_v^\p$, such that
\begin{eqnarray*}
z\cdot\biggl(\Bigl(\det\bigl(\bigoplus_{v\in V}
\psi_v^{\oplus\sigma_v} \bigr)^{\otimes
c}\otimes{\id}_{\Oh_X(m)}\Bigr)\circ
\tau\circ{\Bigl(\bigl(\bigoplus_{v\in V} \psi_v^{\oplus\sigma_v}
\bigr)^{\otimes a}\Bigr)^{\oplus b}}^{-1}\biggr) &=&\tau^\p.
\end{eqnarray*}
A \it family of $V$-split sheaves of type $\ul{P}$ with a decoration of type $(a,b,c,m)$
parameterized by the scheme $S$ \rm
consists of
\begin{itemize}
\item a tuple $(\E_{S,v},v\in V)$ on $S\times X$
      of $S$-flat families $\E_{S,v}$ of torsion free sheaves on $X$ with Hilbert polynomial
      $P_v$, $v\in V$;
\item a line bundle $\EL_S$ on $S$;
\item a homomorphism
$$
\tau_S\colon
\biggl(\Bigl(\bigoplus_{v\in V}\E_{S,v}^{\oplus \sigma_v}\Bigr)^{\otimes a}\biggr)^{\oplus b}
\lra \det\Bigl(\bigoplus_{v\in V}\E_{S,v}^{\oplus \sigma_v}\Bigr)^{\otimes c}
\otimes\pi_S^*\EL_S\otimes\pi_X^*\Oh_X(m).
$$
\end{itemize}
Two such families will be called \it equivalent\rm, if there are isomorphisms
$\psi_{S,v}\colon  \E_{S,v}\lra \E_{S,v}^\p$, $v\in V$, and $z_S\colon  \EL_S\lra \EL_S^\p$,
such that
\begin{eqnarray*}
\Bigl(\det\bigl(\bigoplus_{v\in V}
\psi_{S,v}^{\oplus\sigma_v} \bigr)^{\otimes
c}\otimes\pi_S^*(z_S)\otimes{\id}_{\pi_X^*\Oh_X(m)}\Bigr)\circ
\tau\circ{\Bigl(\bigl(\bigoplus_{v\in V}
\psi_{S,v}^{\oplus\sigma_v})^{\otimes a} \Bigr)^{\oplus
b}}^{-1} &=& \tau^\p.
\end{eqnarray*}
\subsubsection*{The semistability condition}
Let $(\E_v,v\in V)$ be a $V$-split sheaf. Then --- as agreed upon
before --- a \it weighted filtration of $(\E_v,v\in V)$ \rm is a
pair $(\E_\bullet,\ul{\alpha})$, composed of a filtration
$$
\E_\bullet\colon\qquad 0\subsetneq (\E_1^v,v\in V)\subsetneq \cdots\subsetneq (\E_s^v,v\in V)
\subsetneq (\E_v,v\in V)
$$
of $(\E_v,v\in V)$ by non-trivial, proper, and saturated $V$-split subsheaves
and a tuple $\ul{\alpha}=(\alpha_1,...,\alpha_s)$ of positive rational numbers.
Recall that
$$
M_{\ul{\sigma},\ul{\chi}}(\E_\bullet,\ul{\alpha})=\sum_{j=1}^s\alpha_j\cdot
\Bigl(P_{\ul{\sigma},\ul{\chi}}(\E_v,v\in V)\cdot{\rk}_{\ul{\sigma}}(\E_j^v,v\in V)
-
P_{\ul{\sigma},\ul{\chi}}(\E_j^v,v\in V)\cdot {\rk}_{\ul{\sigma}}(\E_v,v\in V)\Bigr).
$$
The number
$$
\mu\bigl(\E_\bullet,\ul{\alpha};\tau\bigr)
$$
is defined by the formula in Remark~\ref{explain3}, ii).
\par
Now, we call a $V$-split sheaf $(\E_v,v\in V;\tau)$
with a decoration of type $(a,b,c,m)$ \it (semi)\allowbreak sta\-ble \rm (or more precisely
\it $(\ul{\sigma},\ul{\eta},\delta)$-(semi)stable\rm),
if for every weighted filtration $(\E_\bullet,\ul{\alpha})$ of $(\E_v,v\in V)$
$$
M_{\ul{\sigma},\ul{\chi}}\bigl(\E_\bullet,\ul{\alpha}\bigr)+\delta\cdot
\mu\bigl(\E_\bullet,\ul{\alpha};\tau\bigr)
(\succeq) 0.
$$
\subsubsection*{The first main result}
We define the functors
\begin{eqnarray*}
\ul{\rm M}^{(\ul{\sigma},\ul{\eta},\delta)-\rm (s)s}_{\ul{P}/a/b/c/m}\colon
\q \ul{\rm Sch}_\C &\lra &\ul{\rm Sets}
\\
S &\lma &\left\{\,\hbox{Equivalence classes of families of (semi)stable $V$-split\,}
\atop \,\hbox{sheaves of type $\ul{P}$ with a decoration of type $(a,b,c,m)$}\,\right\}.
\end{eqnarray*}
\begin{Thm}
\label{main1} {\rm i)} There exist a projective scheme ${\cal
M}:={\cal M}^{(\ul{\sigma},\ul{\eta},\delta)-\rm
ss}_{\ul{P}/a/b/c/m}$ and a natural transformation
$\theta\colon  \ul{\rm M}^{(\ul{\sigma},\ul{\eta},\delta)-\rm
ss}_{\ul{P}/a/b/c/m} \lra h_{\cal M}$, such that for any other
scheme ${\cal M}^\p$ and any other natural transformation
$\theta^\p\colon \ul{\rm M}^{(\ul{\sigma},\ul{\eta},\delta)-\rm
ss}_{\ul{P}/a/b/c/m}\lra h_{\cal M^\p}$, there is a unique
morphism $\zeta\colon  \cal M\lra \cal M^\p$ with
$\theta^\p=h_\zeta\circ\theta$.
\par
{\rm ii)} The space ${\cal M}$ contains an open subscheme ${\cal M}^{\rm s}$
which is a coarse moduli scheme for the functor
$\ul{\rm M}^{(\ul{\sigma},\ul{\eta},\delta)-\rm s}_{\ul{P}/a/b/c/m}$.
\end{Thm}
\subsection{Applications to quiver problems}
\label{appli}
Let $Q=(V,A,t,h)$ be a quiver with vertices
$V=\{\,v_1,...,v_t\,\}$, arrows $A=\{\,a_1,...,a_u\,\}$, the tail
map $t\colon A\lra V$, and the head map $h\colon A\lra V$. We
assume that no multiple arrows occur. Fix a tuple of coherent
sheaves $\ul{\G}=(\G_a,a\in A)$. An \it (augmented) representation
of $Q$ of type $(\ul{P},\ul{\G})$ \rm is a tuple $(\E_v,v\in
V;f_a, a\in A;\eps)$, consisting of
\begin{itemize}
\item a $V$-split sheaf $(\E_v,v\in V)$ of type $\ul{P}$,
\item homomorphisms $f_a\colon \G_{a}\otimes\E_{t(a)}\lra \E_{h(a)}$, $a\in A$,
\item a complex number $\eps$,
\end{itemize}
such that either $\eps\neq 0$ or one of the $f_a$, $a\in A$, is
non trivial. For simplicity, we will often drop the term
``augmented" in the following.
Two representations $(\E_v,v\in V;f_a, a\in A;\eps)$ and
$(\E^\p_v,v\in V;f^\p_a, a\in A;\eps)$ are called \it equivalent\rm, if there
are isomorphisms $\psi_v\colon  \E_v\lra\E_v^\p$, $v\in V$, and $z\in\C^*$,
such that
$$
z\cdot \Bigl(\psi_{h(a)}
\circ f_a\circ \bigl({\id}_{\G_{a}}\otimes\psi_{t(a)}\bigr)^{-1}\Bigr)= f_a^\p,
\q a\in A,\qquad\hbox{and}\q z\cdot\eps=\eps^\p.
$$
A \it family of representations of $Q$ of type $(\ul{P},\ul{\G})$ parameterized
by $S$ \rm is a tuple
$(\E_{S,v},v\in V;f_{S,a}, a\in A; \EL_S, \eps_S)$ \rm which consists of
\begin{itemize}
\item $S$-flat families $\E_{S,v}$ on $S\times X$ of torsion free sheaves on $X$
      with Hilbert polynomial $P_v$, $v\in V$;
\item a line bundle $\EL_S$ on $S$;
\item a section $\eps_S\in H^0(S,\EL_S)$;
\item homomorphisms
      $
      f_{S,a}\colon  \pi_X^*\G_{a}\otimes\E_{S,t(a)}\lra \E_{S,h(a)}\otimes\pi_S^*\EL_S.
      $
\end{itemize}
An \it equivalence of the families \rm $(\E^1_{S,v},v\in
V;f^1_{S,a}, a\in A; \EL^1_S, \eps^1_S)$ and $(\E^2_{S,v},v\in
V;f^2_{S,a}, a\in A; \EL^2_S, \eps^2_S)$ consists of an isomorphism
$z_S\colon \EL_S^1\lra \EL_S^2$ and isomorphisms $\psi_v\colon
\E_{S,v}^1\lra \E_{S,v}^2$, $v\in V$, such that
$$
f_{S,a}^2=\bigl(\psi_{h(a)}\otimes \pi_S^*(z_S)\bigr)\circ
f_{S,a}^1\circ \bigl({\id}_{\pi_X^*\G_a}\otimes
\psi_{t(a)}\bigr)^{-1},\q a\in A,\qquad z_S\circ\eps_S^1=\eps_S^2.
$$
\subsubsection*{Associated decorations}
Fix the parameters $\ul{\sigma}=(\sigma_v\in\Z_{>0}, v\in V)$.
There are an $m\ge 0$ and $b>0$, such that we have surjections
$\nu_a\colon \Oh_C(-m)^{\oplus b}\lra\G_a$ for all $a\in A$, and an embedding
$v_0\colon  \Oh_X\subset\Oh_X(m)$.
Set $M:=\bigoplus_{v\in V}\C^{\sigma_v\cdot r_v}$.
Next, decompose the $\GL(M)$-module
$Z:=\bigl(M^{\otimes\dim M}\otimes (\bigwedge^{\dim M} M)^{\otimes -1}\bigr)^{\oplus (b+1)}$
as $\Hom(M,M)^{\oplus b}\oplus\C\oplus Z^\p$.
\par
Let $(\E_{S,v},v\in V;f_{S,a}, a\in A; \EL_S, \eps)$ be a family of representations
on $S\times X$ and $\iota\colon  U\subset S\times X$ the maximal open subset
over which all the $\E_{S,v}$, $v\in V$, are locally free.
Then, the restrictions of the $f_{S,a}^{\oplus\sigma_{t(a)}\cdot\sigma_{h(a)}}$
to $U$ together with
the pullback of $\eps_S$ to $U$ may be interpreted as
a homomorphism
$$
\tau^\p\colon\ul{\rm Hom}\Bigl(
\bigoplus_{v\in V}\E_{S,v|U}^{\oplus \sigma_v},
\bigoplus_{v\in V}\E_{S,v|U}^{\oplus \sigma_v}\Bigr)^{\oplus b}
\oplus\Oh_U\lra \iota^*\bigl(\pi_S^*\EL_S
\otimes\pi_X^*\Oh_X(m)\bigr).
$$
The splitting of $Z$ yields a natural projection
$$
\det
\Bigl(\E_{S,v|U}^{\oplus \sigma_v}\Bigr)^{\otimes -1}\otimes \biggl(
\Bigl(\bigoplus_{v\in V}\E_{S,v|U}^{\oplus\sigma_v}\Bigr)^{\otimes \sum_{v\in V}\sigma_v\cdot
r_v}\biggr)^{\oplus(b+1)}
\lra
\ul{\rm Hom}\Bigl(
\bigoplus_{v\in V}\E_{S,v|U}^{\oplus \sigma_v},
\bigoplus_{v\in V}\E_{S,v|U}^{\oplus \sigma_v}\Bigr)^{\oplus b}\oplus\Oh_U
,
$$
so that we get
$$
\tau^{\p\p}\colon
\biggl(
\Bigl(\bigoplus_{v\in V}\E_{S,v|U}^{\oplus\sigma_v}\Bigr)^{\otimes \sum_{v\in V}\sigma_v\cdot
r_v}\biggr)^{\oplus (b+1)}\lra \iota^*\biggl(
\det
\Bigl(\E_{S,v}^{\oplus \sigma_v}\Bigr)\otimes\pi_S^*\EL_S
\otimes\pi_X^*\Oh_X(m)\biggr).
$$
We finally define
\begin{eqnarray*}
&\tau_S\colon& \biggl( \Bigl(\bigoplus_{v\in
V}\E_{S,v}^{\oplus\sigma_v}\Bigr)^{\otimes \sum_{v\in
V}\sigma_v\cdot r_v}\biggr)^{\oplus (b+1)} \q\lra\q \iota_*
\biggl( \Bigl(\bigoplus_{v\in
V}\E_{S,v|U}^{\oplus\sigma_v}\Bigr)^{\otimes \sum_{v\in
V}\sigma_v\cdot r_v}\biggr)^{\oplus (b+1)}\q\lra
\\
&&\hskip -10pt \stackrel{\iota_*(\tau^{\p\p})}{\lra}\q
\iota_*\iota^*\biggl(
\det
\Bigl(\E_{S,v}^{\oplus \sigma_v}\Bigr)\otimes\pi_S^*\EL_S
\otimes\pi_X^*\Oh_X(m)\biggr)\stackrel{\rm Prop.\ \ref{Mar}}{=}
\det
\Bigl(\E_{S,v}^{\oplus \sigma_v}\Bigr)\otimes\pi_S^*\EL_S
\otimes\pi_X^*\Oh_X(m).
\end{eqnarray*}
We call $(\E_{S,v},v\in V;\tau_S)$ the \it associated family of
$V$-split sheaves with a decoration of type $(s,b+1,1,m)$, $s:=\sum_{v\in V}\sigma_v\cdot
r_v$\rm. When $S$ is just a point, we simply speak of the \it associated decoration \rm $\tau$.
\subsubsection*{The semistability condition}
Fix the same data as before.
We call a representation $(\E_v,v\in V; f_a, a\in A; \eps)$ \it
$(\ul{\sigma},\ul{\eta},\delta)$-(semi)sta\-ble\rm, if the
$V$-split sheaf $(\E_v,v\in V; \tau)$ with the associated decoration
is $(\ul{\sigma},\ul{\eta},\delta)$-(semi)stable.
\subsubsection*{The second main result}
Define the functors
\begin{eqnarray*}
\ul{\rm R}(Q)^{(\ul{\sigma},\ul{\eta},\delta)-\rm (s)s}_{\ul{P}/\ul{\G}}\colon
\q \ul{\rm Sch}_\C &\lra &\ul{\rm Sets}
\\
S &\lma &\left\{\,\hbox{Equivalence classes of families of $(\ul{\sigma},\ul{\eta},\delta)$-\,}
\atop\,\hbox{(semi)stable representations of $Q$ of type $(\ul{P},\ul{\G})$}\,\right\}.
\end{eqnarray*}
\begin{Thm}
\label{main2}
{\rm i)} There exist a projective scheme
${\cal R}:={\cal R}(Q)^{(\ul{\sigma},\ul{\eta},\delta)-\rm ss}_{\ul{P}/\ul{\G}}$ and
a natural transformation
$\theta\colon  \ul{\rm R}(Q)^{(\ul{\sigma},\ul{\eta},\delta)-\rm ss}_{\ul{P}/\ul{\G}}
\lra h_{\cal R}$, such that for any other scheme ${\cal R}^\p$ and any other
natural transformation $\theta^\p\colon
\ul{\rm R}(Q)^{(\ul{\sigma},\ul{\eta},\delta)-\rm ss}_{\ul{P}/\ul{\G}}\lra h_{\cal R^\p}$,
there exists a unique morphism $\zeta\colon  \cal R\lra \cal R^\p$
with $\theta^\p=h_\zeta\circ\theta$.
\par
{\rm ii)} The scheme ${\cal R}$ contains an open subscheme ${\cal R}^{\rm s}$
which is a coarse moduli scheme for the functor
$\ul{\rm R}(Q)^{(\ul{\sigma},\ul{\eta},\delta)-\rm s}_{\ul{P}/\ul{\G}}$.
\end{Thm}
\subsection{Behaviour for large $\delta$}
\label{large}
Intuitively, one would like to have a semistability concept for representations
of quivers which poses conditions on subrepresentations only.
However, as illustrated in \cite{Sch} for the example of the quiver consisting of one
vertex and an arrow, connecting the vertex to itself, this property
cannot be expected for general $\delta$. However, as in the case of the aforementioned
quiver, for large $\delta$, the semistability concept will stabilize to one
which is a condition on subrepresentations only. Another nice feature is that,
for large $\delta$, one has a generalized Hitchin map.
Our first result is
\begin{Thm}
\label{stab}
Fix the data $\ul{\sigma}$, $\ul{\eta}$, $\ul{P}$, and $\ul{\G}$.
Let $\delta$ be a positive polynomial of degree exactly $\dim X-1$.
Then, there exists a natural number $n_\infty$, such that for all $n\ge n_\infty$,
the following two conditions on a
representation $(\E_v,v\in V; f_a, a\in A; \eps)$ of $Q$ of type $(\ul{P},\ul{\G})$
are equivalent
\begin{itemize}
\item[\rm 1.] The representation
$(\E_v,v\in V; f_a, a\in A; \eps)$ is $(\ul{\sigma},\ul{\eta}/n,n\cdot \delta)$-(semi)stable,
$\ul{\eta}/n:=(\eta_v/n,v\in V)$.
\item[\rm 2.]
\begin{itemize}
\item[\rm a.]
For any non-trivial, proper
subrepresentation $(\F_v,v\in V)$ (i.e., $V$-split subsheaf, such that
$f_{a}(\G_{a}\otimes\F_{t(a)})\subset \F_{h(a)}$ for all $a\in A$)
one has
$$
\frac{P_{\ul{\sigma},\ul{\chi}}(\F_v,v\in
V)}{\rk_{\ul{\sigma}}(\F_v,v\in V)} \q (\preceq) \q
\frac{P_{\ul{\sigma},\ul{\chi}}(\E_v,v\in
V)}{\rk_{\ul{\sigma}}(\E_v,v\in V)},
$$
and
\item[\rm b.]
either $\eps\neq 0$, or the restriction of $(\E_v,v\in V; f_a, a\in A; \eps)$ to a general
point $x\in X$ in the open subset where all the $\E_v$, $v\in V$, are locally free is
semistable.
\end{itemize}
\end{itemize}
\end{Thm}
We call a representation $(\E_v,v\in V; f_a, a\in A; \eps)$ of $Q$ of type $(\ul{P},
\ul{\G})$ \it $(\ul{\sigma},\ul{\chi})$-(semi)stable\rm, if it satisfies Condition 2
in Theorem~\ref{stab} with $(\preceq)$. Note that, for representations of the form
$(\E_v,v\in V; f_a, a\in A; \eps=1)$, this is exactly the ``Gieseker-analogue"
of the semistability
definition given by \'Alvarez-C\'onsul and Garc\'\i a-Prada \cite{AG}.
Observe that for $(\ul{\sigma},\ul{\chi})$-semistable representations, one has the concepts of
a \it Jordan-H\"older filtration\rm, \it the associated graded object\rm,
and \it S-equivalence\rm. Therefore, one can also speak of
\it $(\ul{\sigma},\ul{\chi})$-polystable
representations\rm.
\subsubsection*{Invariants of quivers and the generalized Hitchin map}
Recall that we may find $b>0$ and $m\ge 0$, such that there are surjections
$\nu_a\colon \Oh_X(-m)^{\oplus b}\lra\G_a$ for all $a\in A$ as well as an embedding
$\nu_0\colon \Oh_X\lra\Oh_X(m)$. Therefore, we look now at the quiver $Q^b:=
(V,A^b:=A\times\{\,1,...,b\,\},t,h)$ where the tail and head maps are given by the
projection onto the first factor followed by the tail and head map of $Q$.
In other words, any arrow in $Q$ is replaced by $b$ copies of the same arrow.
We choose $\ul{r}=(r_v,v\in V)$ as the dimension vector. The variety of
representations with this dimension vector is thus
$$
{\rm Rep}_{\ul{r}}(Q^b):=\bigoplus_{a\in
A^b}\Hom\bigl(\C^{r_{t(a)}},\C^{r_{h(a)}}\bigr).
$$
The variety ${\rm Rep}_{\ul{r}}(Q^b)$ comes with an action of the group
$\mathop{\sf X}_{v\in V}\GL(r_v,\C)$. By the work of LeBruyn and Procesi \cite{LP},
one knows explicit generators for the ring of $(\mathop{\sf X}_{v\in V}\GL(r_v,\C))$-invariant
regular functions on ${\rm Rep}_{\ul{r}}(Q^b)$. To state the result, let
${\frak o}=(a_1,...,a_t)$ be an \it oriented cycle\rm, i.e., a sequence
of arrows $a_1,...,a_t$, such that $h(a_{i})=t(a_{i-1})$, $i=2,...,l$,
and $h(a_1)=t(a_t)=:v({\frak o})$. We call
$l$ the \it length of the cycle\rm. For any such cycle and any point
$x\in {\rm Rep}_{\ul{r}}(Q^b)$, we get an endomorphism
$f_{x,\frak o}\colon \C^{r_{v(\frak o)}}\lra\C^{r_{v(\frak o)}}$. We then define
\begin{eqnarray*}
t_{\frak o}\colon {\rm Rep}_{\ul{r}}(Q^b) &\lra&\C
\\
                      x &\lma& {\rm Trace}(f_{x,\frak o}).
\end{eqnarray*}
The function $t_{\frak o}$ is obviously invariant under the
$(\mathop{\sf X}_{v\in V}\GL(r_v,\C))$-action. The result of \cite{LP} states
that the invariants of the form $t_{\frak o}$, $\frak o$ an oriented cycle,
generate the ring of invariants
${\cal R}_{\ul{r}}(Q^b):=\C[{\rm Rep}_{\ul{r}}(Q^b)]^{\mathop{\sf X}_{v\in V}\GL(r_v,\C)}$. 
Moreover, one may restrict to oriented cycles of length
at most $\bigl(\sum_{v\in V} r_v\bigr)^2+1$.
We also look at the affine variety
$$
{\rm Rep}^\eps_{\ul{r}}(Q^b):= {\rm Rep}_{\ul{r}}(Q^b)\oplus\C.
$$
This is also a $(\mathop{\sf X}_{v\in V}\GL(r_v,\C))$-variety,
the action on $\C$ being trivial. Denote by $t_0$ the projection onto the second
factor. This is a $(\mathop{\sf X}_{v\in V}\GL(r_v,\C))$-invariant function,
and the above result implies
\begin{eqnarray*}
{\cal R}^\eps_{\ul{r}}(Q^b)
&:=&
\C\bigl[{\rm Rep}^\eps_{\ul{r}}(Q^b)\bigr]^{\mathop{\sf X}_{v\in V}\GL(r_v,\C)}
\\
&=&
\C\Bigl[t_0;\ t_{\frak o}, {\frak o} \hbox{ an oriented cycle of length} \le
\bigl(\sum_{v\in V} r_v\bigr)^2+1\Bigr].
\end{eqnarray*}
Next, assign to $t_0$ the degree one and to $t_{\frak o}$ the degree
length of ${\frak o}$, $\frak o$ an oriented cycle.
Then, ${\cal R}^\eps_{\ul{r}}(Q^b)$ is a graded ring, and
${\rm Proj}({\cal R}^\eps_{\ul{r}}(Q^b))$ identifies with the
$(\C^*\times \mathop{\sf X}_{v\in V}\GL(r_v,\C))$-quotient of
${\rm Rep}^\eps_{\ul{r}}(Q^b)$ where $\C^*$ acts by scalar multiplication.
We may choose a degree $d$ such that
the subring ${\cal R}^{(d)}\subset {\cal R}^\eps_{\ul{r}}(Q^b)$ of elements the degree
of which is a multiple of $d$ is generated by elements of degree $d$,
say, ${\frak i}_0,...,{\frak i}_q$ (see \cite{Red}, III.8, Lemma).  This yields an embedding
$$
{\rm Proj}({\cal R}^\eps_{\ul{r}}(Q^b))
\q\hookrightarrow\q
\P_{q}.
$$
Now, we return to the setting of representations of $Q$ of type $(\ul{P},\ul{\G})$
where we fix $b$, $m$, $\nu_a$, $a\in A$, and $\nu_0$ as before.
Set ${\Bbb H}(Q,\ul{P},\ul{\G}):=\P\bigl(H^0(\Oh_X(d\cdot m))^{\oplus q}\bigr)$.
\par
Let $S$ be a scheme, and $(\E_{S,v},v\in V; f_{S,a}, a\in A;\EL_S, \eps)$
a family of $(\ul{\sigma},\ul{\chi})$-semistable representations of $Q$ of
type $(\ul{P},\ul{\G})$ parameterized by $S$. Denote by $\iota\colon U\subset S\times X$ the
maximal open subset where all the $\E_{S,v}$, $v\in V$, are locally free.
To the invariant $t_0$ corresponds the homomorphism
$$
t^{\p\p}_0\colon \Oh_{S\times X}\stackrel{\pi_S^*\eps_S}{\lra}\pi_S^*\EL_S
\stackrel{{\id}_{\pi_S^*\EL_S}\otimes\pi_X^*\nu_0}{\lra} \pi_S^*\EL_S\otimes\pi_X^*\Oh_X(m).
$$
Furthermore, using the quotients $\nu_a$, $a\in A$, for any oriented cycle
${\frak o}$ of length $l$, we get a homomorphism
$$
f^\p_{\frak o}\colon \E_{S,v(\frak o)}\lra \E_{S,v(\frak o)}\otimes
\bigl(\pi_S^*\EL_S\otimes
\pi_X^*\Oh_X(m)\bigr)^{\otimes l}.
$$
If we restrict $f^\p_{\frak o}$ to $U$ and take traces, we obtain a section
$$
t^\p_{\frak o}\colon \Oh_U\lra \iota^*\bigl(\pi_S^*\EL_S\otimes
\pi_X^*\Oh_X(m)\bigr)^{\otimes l}.
$$
By Proposition \ref{Mar}, this extends to
$$
t^{\p\p}_{\frak o}\colon
\Oh_{S\times X}\lra \bigl(\pi_S^*\EL_S\otimes
\pi_X^*\Oh_X(m)\bigr)^{\otimes l}.
$$
Therefore, any invariant ${\frak i}_j$, $j=0,...,q$, provides a section
$$
i_j:=i_j\bigl(\E_{S,v},v\in V; f_{S,a}, a\in A;\EL_S, \eps\bigr)\colon \Oh_{S\times X}
\lra \bigl(\pi_S^*\EL_S\otimes
\pi_X^*\Oh_X(m)\bigr)^{\otimes d},\q j=0,...,q.
$$
Condition 2.b now grants that one of the homomorphisms $i_j$, $j=0,...,q$, will
be non-zero.
Hence, we get a morphism
$$
{\rm Hit}\bigl(\E_{S,v},v\in V; f_{S,a}, a\in A;\EL_S, \eps\bigr)
\colon S\lra {\Bbb H}:={\Bbb H}(Q,\ul{P},\ul{\G})
$$
with ${\rm Hit}^*\Oh_{\Bbb H}(1)=\EL_S^{\otimes d}$.
\subsubsection*{The third main result}
This time, we look at the functors
\begin{eqnarray*}
\ul{\rm R}(Q)^{(\ul{\sigma},\ul{\chi})-\rm (s)s}_{\ul{P}/\ul{\G}}\colon
\q \ul{\rm Sch}_\C &\lra &\ul{\rm Sets}
\\
S &\lma &\left\{\,\hbox{Equivalence classes of families of $(\ul{\sigma},\ul{\chi})$-\,}
\atop\,\hbox{(semi)stable representations of $Q$ of type $(\ul{P},\ul{\G})$}\,\right\}.
\end{eqnarray*}
\begin{Thm}
\label{main3}
{\rm i)} There exist a projective scheme
${\cal R}:={\cal R}(Q)^{(\ul{\sigma},\ul{\chi})-\rm ss}_{\ul{P}/\ul{\G}}$ and
a natural transformation
$\theta\colon  \ul{\rm R}(Q)^{(\ul{\sigma},\ul{\chi})-\rm ss}_{\ul{P}/\ul{\G}}
\lra h_{\cal R}$, such that for any other scheme ${\cal R}^\p$ and any other
natural transformation $\theta^\p\colon
\ul{\rm R}(Q)^{(\ul{\sigma},\ul{\chi})-\rm ss}_{\ul{P}/\ul{\G}}\lra h_{\cal R^\p}$,
there exists a unique morphism $\zeta\colon  \cal R\lra \cal R^\p$
with $\theta^\p=h_\zeta\circ\theta$.
\par
{\rm ii)} The scheme ${\cal R}$ contains an open subscheme ${\cal R}^{\rm s}$
which is a coarse moduli scheme for the functor
$\ul{\rm R}(Q)^{(\ul{\sigma},\ul{\chi})-\rm s}_{\ul{P}/\ul{\G}}$.
\par
{\rm iii)} The closed points of ${\cal R}$ are in bijection to the set
of S-equivalence classes of $(\ul{\sigma},\ul{\chi})$-semistable representations
of $Q$ of type $(\ul{P},\ul{Q})$, or, equivalently, to the set of isomorphism
classes of $(\ul{\sigma},\ul{\chi})$-polystable representations
of $Q$ of type $(\ul{P},\ul{Q})$.
\par
{\rm iv)} There is a {\rm generalized Hitchin morphism}
$$
{\rm Hit}(Q,\ul{P},\ul{\G})\colon
\cal R\lra {\Bbb H}(Q,\ul{P},\ul{\G}).
$$
\end{Thm}
Note that the theorem in the introduction now follows by taking
${\cal D}$ and ${\Bbb D}$ as the open subscheme $\eps=1$ in
${\cal R}$ and ${\Bbb H}$, respectively (cf.\ Remark~\ref{explain3}).
\section{The proofs}
\subsection{Proof of Theorem~\ref{main1}}
\label{proofi}
In order to prove Theorem~\ref{main1}, one can copy almost word by word the
proofs in \cite{Sch}, Section~2.3.6, or \cite{GS}. The only point which has to be given special
attention and which is indeed rather tricky is the correct choice of a linearization on the parameter space.
We will, therefore, construct the parameter space, give the linearization
of the respective group action, and show in a sample computation that it is the
correct one.
\subsubsection*{The parameter space}
We denote by ${\frak A}_v$, $v\in V$,
the union of those components of ${\rm Pic}(X)$
which contain line bundles of the form $\det(\E_v)$
where
$(\E_v,v\in V;\tau)$
is a semistable $V$-split sheaf of type
$\ul{P}$ with a decoration of type $(a,b,c,m)$.
We also set ${\frak A}:=\mathop{\sf  X}_{v\in V}{\frak A}_v$.
By the usual boundedness arguments, we can find an $l_0$, such that for all
$l\ge l_0$, all semistable $V$-split sheaves
$(\E_v,v\in V; \tau)$ of type
$\ul{P}$ with a decoration of type $(a,b,c,m)$,
all $v\in V$, all $[\EL]\in {\frak A}_v$, and all
$\n=\bigotimes_{v\in V}\EL_v^{\otimes \sigma_v}$ with $[\EL_v]\in {\frak A}_v$,
$v\in V$,
\begin{itemize}
\item $\E_v(l)$ is globally generated and $H^i(\E_v(l))=0$ for all $i>0$;
\item $\EL(r_v\cdot l)$ is globally generated and $H^i(\EL(r_v\cdot l))=0$ for all $i>0$;
\item $\n^{\otimes c}(a\cdot l)$ is globally generated and $H^i(\n^{\otimes c}(a\cdot l))=0$ for all $i>0$.
\end{itemize}
We fix such an $l$, and set $p_v:=P_v(l)$, $v\in V$, and $p:=\sum_{v\in V}\sigma_v\cdot p_v$.
Moreover, we choose vector spaces $W_v$ of dimension $p_v$ and let
${\frak Q}^0_v$ be the quasi-projective quot scheme parameterizing
quotients $q\colon W_v\otimes\Oh_X(-l)\lra \F$ with $\F$ a torsion free coherent
$\Oh_X$-module with Hilbert polynomial $P_v$ and $H^0(q(l))$ an isomorphism, $v\in V$.
Let ${\frak E}_v$ be the universal quotient on ${\frak Q}_v^0\times X$, $v\in V$,
and
$$
{\frak E}^{\rm total}:=\bigoplus_{v\in V}\pi_{\frak Q^0_v}^*{\frak E}_v^{\oplus\sigma_v}
$$
be the resulting sheaf on $(\mathop{\sf X}_{v\in V}{\frak Q}_v^0)\times X$.
Define $M:=\bigoplus_{v\in V} W^{\oplus\sigma_v}_v$ and
$$
{\frak P}:= \P \Bigl(\bigl((M^{\otimes a})^{\oplus b}\bigr)^\vee\otimes\pi_{
(\mathop{\sf X}_{v\in V}{\frak Q}_v^0)*}
\bigl(\det({\frak E^{\rm total}})^{\otimes c}\otimes\pi_X^*\Oh_X(a\cdot l)\bigr)\Bigr).
$$
This is a projective bundle over $\mathop{\sf X}_{v\in V}{\frak Q}_v^0$,
and the parameter space ${\frak M}$ is constructed in the usual
way as a closed subscheme of ${\frak P}$.
In particular, it is projective over $\mathop{\sf X}_{v\in V}{\frak Q}_v^0$.
Furthermore, ${\frak M}$ comes with an action of the group
$\bigl(\mathop{\sf X}_{v\in V}\GL(W_v)\bigr)/\C^*$, $\C^*$ being diagonally embedded.
We define
\begin{eqnarray*}
\widetilde{G}
&:=&
\bigl(\mathop{\sf X}_{v\in V}\GL(W_v)\bigr)\cap\SL(M)\\
&=&
\Bigl\{\,(h_1,...,h_t)\in\mathop{\sf X}_{v\in V}\GL(W_v)
\,|\,
\det(h_1)^{\sigma_{v_1}}\cdot...\cdot\det(h_t)^{\sigma_{v_t}}=1\,\Bigr\}.
\end{eqnarray*}
The group $\widetilde{G}$ maps with finite kernel onto $\bigl(\mathop{\sf X}_{v\in V}\GL(W_v)\bigr)/\C^*$,
whence we may restrict our attention to the action of $\widetilde{G}$.
\par
The linearization of the above group action will be induced via a Gieseker morphism
to some other scheme. For this, we fix Poincar\'e line bundles ${\cal P}_v$
over ${\frak A}_v\times X$, $v\in V$, and set
$$
{\frak G}_v:=\P\Bigl(\bigl(\bigwedge^{r_v}W_v\bigr)^\vee \otimes
\pi_{\frak A_v *}\bigl({\cal P}_v\otimes\pi_X^*\Oh_X(r_v\cdot l)\bigr)\Bigr).
$$
Choosing ${\cal P}_v$ appropriately, we may assume that $\Oh_{\frak G_v}(1)$ is very ample
for all $v\in V$. On ${\frak A}\times X$, we get the line bundle
$$
{\cal P}:=\bigotimes_{v\in V}\pi_{\frak A_v\times X}^*{\cal P}_v^{\otimes\sigma_v}.
$$
Then, we define
$$
{\frak P}^\p:= \P \Bigl(\bigl((M^{\otimes a})^{\oplus b}\bigr)^\vee\otimes\pi_{\frak A*}
\bigl({\cal P}^{\otimes c}\otimes\pi_X^*\Oh_X(a\cdot l)\bigr)\Bigr)
$$
as a projective bundle over ${\frak A}$.
Again, $\Oh_{\frak P^\p}(1)$ can be assumed to be ample.
We now have a $\widetilde{G}$-equivariant and injective
morphism
$$
\Gamma\colon {\frak M}\lra {\frak P}^\p\times \mathop{\sf X}_{v\in V} {\frak G}_v.
$$
For given $\beta\in\Z_{>0}$, and $\kappa_v\in\Z_{>0}$, $v\in V$, there is a natural
linearization of the $\widetilde{G}$-action on
${\frak P}^\p\times \mathop{\sf X}_{v\in V} {\frak G}_v$ in the ample line bundle
$\Oh(\beta;\kappa_v,v\in V)$. This may be altered by any character of
$\mathop{\sf X}_{v\in V}\GL(W_v)$.
Let $d:=\delta(l)$, $x_v:=-\chi_v(l)/d$, $x_v^\p:=r_v x_v/p_v$, $v\in V$,
$$
\eps:= \frac{p-a\cdot d}{r\cdot d},\q
\eps_v:=\sigma_v-\frac{x_v}{\eps}=\sigma_v+\frac{r\cdot \chi_v(l)}{p-a\cdot d},\q v\in V,
$$
and
$$
x_\nu^{\p\p}:= \eps\cdot\sigma_\nu\cdot
\left(\frac{r}{p}-\frac{r_v}{p_v}\right),\q v\in V.
$$
\begin{Rem}
To be very precise, the quantities $\eps$ and $\eps_v$, $v\in V$, are functions in $l$.
Since $p=P(l)$ is a positive polynomial of degree $\dim X$ and both $\delta$ and $\chi_v$
are polynomials of degree at most $\dim X-1$, it is clear that $\eps$ and $\eps_v$, $v\in V$,
will be positive for $l\gg 0$, i.e., the line bundle in which the action is linearized is
really ample.
\end{Rem}
Now, we choose $\beta\in\Z_{>0}$ and $\kappa_v\in\Z_{>0}$
such that
$$
\frac{\kappa_v}{\beta}= \eps\cdot\eps_v,\q v\in V.
$$
We modify the linearization of the $\widetilde{G}$-action on
$\mathop{\sf X}_{v\in V}{\frak G}_v$ in $\Oh(\kappa_v,v\in V)$ by
a character, such that
${\C^*}^t=\C^*\cdot{\id}_{W_{v_1}}\times\cdots\times\C^*\cdot{\id}_{W_{v_t}}$
acts via a $(z_v,v\in V)\lma {\Pi}_{v\in V} z_v^{p_v\cdot e_v}$
with
$$
e_v:= \beta\cdot (x_v^\p+x_v^{\p\p}),\q v\in V.
$$
Note that this character is just the restriction of the character
$$
(m_1,...,m_t)\lma 
\det(m_1)^{e_{v_1}}\cdot..\cdot\det(m_t)^{e_{v_t}}
$$
of ${\sf X}_{v\in V}\GL(W_v)$ to the center ${\cal Z}$. We work
with the resulting linearization of the $\widetilde{G}$-action on
${\frak P}^\p\times \mathop{\sf X}_{v\in V}{\frak G}_v$ in
$\Oh(\beta;\kappa_v,v\in V)$.
\subsubsection*{A sample computation}
In order to illustrate that our choice of the linearization is accurate, we go through
a part of the calculations which are analogous to those in Section~2.3 of \cite{Sch}.
More precisely, we show the following: Let $m=(q_v\colon W_v\otimes\Oh_X(-l)\lra\E_v, v\in V;
\tau)$ be a point in the parameter space ${\frak M}$, such that
$\Gamma(m)$ is (semi)stable w.r.t.\ the chosen linearization in $\Oh(\beta;\kappa_v,v\in V)$,
then $(\E_v,v\in V;\tau)$ is a (semi)stable $V$-split sheaf with a decoration
of type $(a,b,c,m)$. First, as in \cite{Sch}, one verifies that the (semi)stability
condition has to be checked only for those weighted filtrations $(\E_\bullet,\ul{\alpha})$
which satisfy
$$
\E^v_j(l)\  \hbox{is globally generated and}\  H^i\bigl(\E^v_j(l)\bigr)=0,
\q i>0,\ j=1,...,s,\ v\in V.
$$
For weighted filtrations of that type, we have to prove that
\begin{equation}
\label{checki}
M_{\ul{\sigma},\ul{\chi}}\bigl(\E_\bullet,\ul{\alpha}\bigr)(l)
+\delta(l)\cdot \mu\bigl(\E_\bullet,\ul{\alpha};\tau\bigr) (\ge)  0.
\end{equation}
Define $\ul{\gamma}=(\gamma_1,...,\gamma_{s+1})$ by the conditions
$$
\frac{\gamma_{j+1}-\gamma_j}{p}= \alpha_j,\q j=1,...,s,
$$
and, setting $\E_j^{\rm total}:=\bigoplus_{v\in V} \E_j^{v,\oplus\sigma_v}$, $j=1,...,s+1$,
$$
\sum_{j=1}^{s+1}\gamma_j\cdot \Bigl(h^0\bigl(\E_j^{\rm total}(l)\bigr)
- h^0\bigl(\E_{j-1}^{\rm total}(l)\bigr)\Bigr)= 0.
$$
Then, we obtain a weighted filtration $(\E_\bullet,\ul{\gamma})$ and, thus,
weighted filtrations $(\widehat{\E}^v_\bullet,\ul{\gamma}^v)$ of the $\E_v$, $v\in V$.
Next, we choose bases $\ul{w}^v=(w_1^v,...,w_{p_v}^v)$ of the $W_v$
with
$$
\langle\,w_1^v,...,w^v_{h^0(\widehat{\E}^v_j(l))}\,\rangle= H^0\bigl(\widehat{\E}^v_j(l)\bigr),
\q j=1,...,s_v,\ v\in V,
$$
and set
$$
\widetilde{\ul{\gamma}}^v:=\bigl(\,\underbrace{\gamma^v_1,...,\gamma^v_1}
_{h^0(\widehat{\E}^v_j(l))\times},...,
\underbrace{\gamma^v_{s_v+1},...,\gamma^v_{s_v+1}}_{(p_v-h^0(\widehat{\E}^v_{s_v}(l)))\times}\,\bigr).
$$
This yields the one parameter subgroup
$$
\la:=\bigl(\la(\ul{w}^{v_1},\widetilde{\ul{\gamma}}^{v_1}),...,
\la(\ul{w}^{v_t},\widetilde{\ul{\gamma}}^{v_t})\bigr)
$$
of $\widetilde{G}$.
Now, with $\Gamma(m)=([m^\p];[m_v],v\in V)$,
\begin{eqnarray*}
\frac{\mu\bigl(\la,\Gamma(m)\bigr)}{\beta} &=& \mu\bigl(\la,
[m^\p]\bigr)
\\
&&+ \eps\cdot \underbrace{\biggl(\sum_{v\in
V}\eps_v\cdot\mu\bigl(\la, [m_v]\bigr) -\sum_{v\in
V}\sigma_v\Bigl(\frac{r_v}{p_v}-\frac{r}{p}\Bigr)\cdot\Bigl(\sum_{j=1}^{s_v+1}\gamma^v_j\bigl(h^0(\widehat{\E}^v_j(l))
-h^0(\widehat{\E}^v_{j-1}(l))\bigr)\Bigr)\biggr)} _{=:A}
\\
&&+\underbrace{\sum_{v\in V}
\Bigl(x_v^\p\cdot\sum_{j=1}^{s_v+1}\gamma^v_j\bigl(h^0(\widehat{\E}^v_j(l))
-h^0(\widehat{\E}^v_{j-1}(l))\bigr)\Bigr)}_{=:B}.
\end{eqnarray*}
Observe
\begin{eqnarray*}
\sum_{j=1}^{s_v+1}\gamma^v_j\bigl(h^0(\widehat{\E}^v_j(l))
-h^0(\widehat{\E}^v_{j-1}(l))\bigr)&=&
\gamma^v_{s_v+1}\cdot p_v-\sum_{j=1}^{s_v}  \bigl(\gamma^v_{j+1}-\gamma^v_j\bigr)
\cdot h^0\bigl(\widehat{\E}^v_j(l)\bigr)
\\
&=&
\gamma_{s+1}\cdot p_v-\sum_{j=1}^{s_v}  \bigl(\gamma^v_{j+1}-\gamma^v_j\bigr)
\cdot P\bigl(\widehat{\E}^v_j\bigr)(l).
\end{eqnarray*}
As $\sum_{v\in V} x_v^\p\cdot p_v=0$, we find
$$
B=
-\sum_{v\in V}\Bigl(x_v^\p\cdot
\sum_{j=1}^{s_v}  \bigl(\gamma^v_{j+1}-\gamma^v_j\bigr)
\cdot P\bigl(\widehat{\E}^v_j\bigr)(l)
\Bigr).
$$
Next,
\begin{eqnarray*}
\mu\bigl(\la,[m_v]\bigr)
&=&
\sum_{j=1}^{s_v} \frac{\gamma^v_{j+1}-\gamma^v_j}{p_v}
\cdot\Bigl(p_v\cdot \rk\widehat{\E}^v_j-h^0\bigl(\widehat{\E}^v_j(l)\bigr)\cdot r_v\Bigr)
\\
&=&
\sum_{j=1}^{s_v} \frac{\gamma^v_{j+1}-\gamma^v_j}{p_v}
\cdot\Bigl(p_v\cdot \rk\widehat{\E}^v_j-P(\widehat{\E}^v_j)(l)\cdot r_v\Bigr).
\end{eqnarray*}
Thus,
\begin{eqnarray*}
&&\eps\cdot\eps_v\cdot\mu\bigl(\la,[m_v]\bigr)
-
x_v^\p\cdot
\sum_{j=1}^{s_v}  \bigl(\gamma^v_{j+1}-\gamma^v_j\bigr)
\cdot h^0\bigl(\widehat{\E}^v_j(l)\bigr)
\\
&=&
\sum_{j=1}^{s_v} \frac{\gamma^v_{j+1}-\gamma^v_j}{p_v}
\cdot\biggl(\eps\cdot\eps_v\cdot\Bigl(p_v\cdot \rk\widehat{\E}^v_j-P(\widehat{\E}^v_j)(l)\cdot r_v\Bigr)
-x_v^\p\cdot p_v\cdot P(\widehat{\E}^v_j)(l)\biggr)
\\
&=&
\sum_{j=1}^{s_v} \frac{\gamma^v_{j+1}-\gamma^v_j}{p_v}
\cdot\biggl(\eps\cdot\eps_v\cdot\Bigl(p_v\cdot \rk\widehat{\E}^v_j-P(\widehat{\E}^v_j)(l)\cdot r_v\Bigr)
-x_v\cdot r_v\cdot P(\widehat{\E}^v_j)(l)\biggr)
\\
&=&
\sum_{j=1}^{s_v} \frac{\gamma^v_{j+1}-\gamma^v_j}{p_v}
\cdot\biggl(\eps\cdot\sigma_v\cdot\Bigl(p_v\cdot \rk\widehat{\E}^v_j-P(\widehat{\E}^v_j)(l)\cdot r_v\Bigr)
-x_v\cdot p_v\cdot \rk\widehat{\E}^v_j\biggr)
\\
&=&
\eps\cdot\sigma_v\cdot  \sum_{j=1}^{s_v}\biggl(\frac{\gamma^v_{j+1}-\gamma^v_j}{p_v}
\cdot\Bigl(p_v\cdot \rk\widehat{\E}^v_j-P(\widehat{\E}^v_j)(l)\cdot r_v\Bigr)\biggr)
-\sum_{j=1}^{s_v}\Bigl(x_v\cdot (\gamma^v_{j+1}-\gamma^v_j)\cdot \rk\widehat{\E}^v_j\Bigr).
\end{eqnarray*}
For a given vertex $v_0\in V$ and a given index $j_0\in\{\,1,...,s_v\,\}$, let
$s(v_0,j_0)\le s^\p(v_0,j_0)\in\{\,1,...,s\,\}$ be the minimal and the maximal index among those indices $j$ with
$\E^{v_0}_{j}=\widehat{\E}^{v_0}_{j_0}$. Then, by definition,
$$
p\cdot \sum_{j=s(v_0,j_0)}^{s^\p(v_0,j_0)}\alpha_{j}= \gamma^{v_0}_{j_0+1}
-\gamma^{v_0}_{j_0}.
$$
Hence,
$$
\sum_{v\in V}
\biggl(\sum_{j=1}^{s_v}\Bigl(x_v\cdot (\gamma^v_{j+1}-\gamma^v_j)\cdot \rk\widehat{\E}
^v_j\Bigr)\biggr)
=
p\cdot\sum_{j=1}^s\alpha_j\Bigl(\sum_{v\in V} \bigl(x_v\cdot \rk\E^v_j\bigr)\Bigr).
$$
Using Proposition~\ref{weightii}, we discover that $\eps\cdot A+B$
equals
\begin{eqnarray*}
&&\eps\cdot\sum_{j=1}^s\alpha_j\Bigl(p\cdot\rk\E^{\rm total}_j-
P\bigl(\E^{\rm total}_j\bigr)(l)\cdot r\Bigr)
-p\cdot\sum_{j=1}^s\alpha_j\Bigl(\sum_{v\in V} \bigl(x_v\cdot
\rk\E^v_j\bigr)\Bigr)
\\
&=&
\sum_{j=1}^s\alpha_j\cdot\biggl(\frac{p^2\rk\E^{\rm total}_j}{r\cdot d}
-\frac{p\cdot a\cdot \rk\E^{\rm total}_j}{r}-\frac{p\cdot P\bigl(
\E^{\rm total}_j\bigr)(l)}{d}
+a\cdot P\bigl(\E^{\rm total}_j\bigr)(l)\biggr)
\\
&&
-p\cdot\sum_{j=1}^s\alpha_j\Bigl(\sum_{v\in V} \bigl(x_v\cdot \rk\E^v_j\bigr)\Bigr).
\end{eqnarray*}
In order to conclude, we have to compute $\mu(\la,[m^\p])$.
Under the identification of $M$ with the space $H^0(\E^{\rm total}(l))$, we define
$$
{\rm gr}_j(M)=H^0\Bigl(\bigl(\E^{\rm total}_j/\E^{\rm total}_{j-1}\bigr)(l)\Bigr),
\q j=1,...,s+1.
$$
The basis $\ul{m}$
of $M$ induced by the bases $\ul{w}^v$ for the $W_v$, $v\in V$, yields a natural
isomorphism
$$
M\cong\bigoplus_{j=1}^{s+1}{\rm gr}_j(M).
$$ For an index tuple
$\ul{\iota}\in J^a:=\{\,1,...,s+1\,\}^{\times a}$, we define
$M_{\ul{\iota}}:={\rm gr}_{\iota_1}(M)\otimes\cdots\otimes{\rm gr}_{\iota_a}(M)$, and
for $k\in\{\,1,...,b\,\}$, we let $M^k_{\ul{\iota}}$ be $M_{\ul{\iota}}$ embedded
into the $k$-th copy of $M^{\otimes a}$ in ${M^{\otimes a}}^{\oplus b}$.
If we denote $P(\E^{\rm total}_j(l))=h^0(\E^{\rm total}_j(l))$ by $m_j$, $j=1,...,s$, then $\la
=\sum_{j=1}^s\alpha_j\la(\ul{m},\ul{\gamma}_p^{(m_j)})$ as a one parameter subgroup
of $\SL(M)$.
Therefore,
$$
\mu(\la,[m^\p])
=
-\min\Bigl\{\,
\sum_{j=1}^s\alpha_j\bigl(a\cdot m_j-\nu_j(\ul{\iota})\cdot p\bigr)
\,|\,
k\in\{\,1,...,b\,\},\ \ul{\iota}\in J^a\colon M^k_{\ul{\iota}}
\not\subset\ker(m^\p)
\,\Bigr\}.
$$
Here,
$$
\nu_j(\ul{\iota})= \#\bigl\{\,\iota_i\le j\,|\, \ul{\iota}=(\iota_1,...,
\iota_a),\ i=1,...,a\,\bigr\}.
$$
Let $\iota_0\in J^a$ be an index which realizes the precise value of $\mu(\la,[m^\p])$.
Then, altogether, we find
\begin{eqnarray*}
&& \sum_{j=1}^s\alpha_j\biggl(\frac{p^2\rk\E^{\rm total}_j}{r\cdot d}
-\frac{p\cdot a\cdot \rk\E^{\rm total}_j}{r}-\frac{p\cdot P\bigl(
\E^{\rm total}_j(l)\bigr)}{d}\biggr)\q+\hskip 4cm
\\
&&\hskip 4cm+\q p\cdot \sum_{j=1}^s \alpha_j\cdot\nu_j(\ul{\iota_0})
-p\cdot\sum_{j=1}^s\alpha_j\Bigl(\sum_{v\in V}\bigl(x_v\cdot \rk\E^v_j\bigr)\Bigr)
\end{eqnarray*}
as the value for $\mu(\la,\Gamma(m))/\beta$. We multiply this by $r\cdot d/p$ and get
\begin{eqnarray*}
&&\sum_{j=1}^s\alpha_j\cdot\biggl(p\cdot \rk\E^{\rm total}_j
-r\cdot P\bigl(\E^{\rm total}_j(l)\bigr)\biggr)
+d\cdot
\Bigl(\sum_{j=1}^s \alpha_j\cdot\bigl(\nu_j(\ul{\iota_0})\cdot r-a\cdot \rk\E^{\rm total}_j
\bigr)\Bigr)+\hskip 1cm
\\
&&
\hskip 7cm+r\cdot\sum_{j=1}^s\alpha_j\Bigl(\sum_{v\in V}\bigl(\chi_v(l)\cdot \rk\E^v_j\bigr)\Bigr).
\end{eqnarray*}
As in \cite{Sch}, one verifies that
$$
\mu\bigl(\E_\bullet,\ul{\alpha};\tau\bigr)
=
\sum_{j=1}^s \alpha_j\cdot\bigl(\nu_j(\ul{\iota_0})\cdot r-a\cdot \rk\E^{\rm total}_j
\bigr),
$$
so that, by Remark~\ref{explain2}, $\mu(\la,\Gamma(m))\ (\ge)\ 0$
implies Inequality~(\ref{checki}).\qed
\subsection{Proof of Theorem~\ref{main2}}
We use the same set up and the same notation as in the Section~\ref{appli}
``Associated decorations" and, w.r.t.\ the corresponding parameters,
at the beginning of Section~\ref{proofi}.
This time, we set $N:=\bigoplus_{v\in V}W_v$.
The space
$$
{\frak P}^{\p\p}
:=
\P\biggl(N^\vee\otimes \pi_{(\mathop{\sf X}_{v\in V}{\frak Q}_v^0)*}
\Bigl(\bigoplus_{v\in V}\pi_{{\frak Q}_v^0}^*{\frak E}_v\otimes
\pi_X^*\Oh_X(l+m)\Bigr)^{\oplus b}
\oplus H^0(X,\Oh_X(m))\otimes\Oh_{\mathop{\sf X}_{v\in V}{\frak Q}_v^0}\biggr)
$$
is a projective bundle over $\mathop{\sf X}_{v\in V}{\frak Q}_v^0$.
Denote by ${\frak E}_{\frak P^{\p\p},v}$ the pullback of ${\frak E}_v$ to
${\frak P}^{\p\p}\times X$. On ${\frak P}^{\p\p}\times X$,
there are the tautological homomorphisms
$$
\phi^{\p\p}
\colon
N\otimes\Oh_{\frak P^{\p\p}\times X}
\lra
\Bigl(\bigoplus_{v\in V}{\frak E}_{\frak P^{\p\p},v}\otimes
\pi_X^*\Oh_X(l+m)\Bigr)^{\oplus b}\otimes\pi_{\frak P^{\p\p}}^*{\Oh}_{\frak P^{\p\p}}(1)
$$
and
$$
\eps^{\p\p}\colon
\Oh_{\frak P^{\p\p}\times X}
\lra
\pi_{\frak P^{\p\p}}^*{\Oh}_{\frak P^{\p\p}}(1)\otimes\pi_X^*\Oh_X(m).
$$
First, we define ${\frak R}^\p$ as the closed subscheme where
$\phi^{\p\p}$ factorizes over the quotient $\bigoplus_{v\in V}{\frak E}_{\frak P^{\p\p},v}
\otimes\pi_X^*\Oh_X(l)$.
Then, $\phi^{\p\p}_{|{\frak R}^\p\times X}$ may be considered as a collection
of homomorphisms
$$
f^\p_{v,v^\p}
\colon
{\frak E}_{\frak R^\p, v} \lra {\frak E}_{\frak R^\p, v^\p}\otimes \pi_X^*\Oh_X(m)^{\oplus b}
\otimes\pi_{\frak R^\p}^*\Oh_{\frak R^\p}(1),\q v,v^\p\in V.
$$
Moreover, we have
$$
\eps^\p\colon
\Oh_{\frak R^\p\times X}\lra \pi^*_{\frak R^\p}\Oh_{\frak R^\p}(1)\otimes\pi_X^*\Oh_X(m).
$$
Now, we can define ${\frak R}$ as a closed subscheme of ${\frak R}^\p$ by the
following conditions
\begin{itemize}
\item The restriction of $f^\p_{v,v^\p}$ to ${\frak R}\times X$ is trivial,
      if $(v,v^\p)\not\in A$;
\item If $(v,v^\p)=a\in A$, then the corresponding homomorphism
$$
f^{\p\p}_{v,v^\p}
\colon
\pi_X^*\Oh_X(-m)^{\oplus b}\otimes{\frak E}_{\frak R^\p, v} \lra {\frak E}_{\frak R^\p, v^\p}\otimes
\pi_{\frak R^\p}^*\Oh_{\frak R^\p}(1)
$$
vanishes
on
$$
\pi_X^*\ker\Bigl(\Oh_X(-m)^{\oplus b}\lra \G_a\Bigr)\otimes {\frak E}_{\frak R^\p, v}.
$$
Note that this is a closed condition, by Proposition \ref{GS}.
\item The restriction of
$$
\Oh_{\frak R^\p\times X}\stackrel{\eps^\p}{\lra}
\pi_X^*\Oh_X(m)\otimes\pi^*_{\frak R^\p}\Oh_{\frak R^\p}(1)
\lra
\pi_X^*\bigl(\Oh_X(m)/\Oh_X\bigr)\otimes\pi^*_{\frak R^\p}\Oh_{\frak R^\p}(1)
$$
to ${\frak R}^\p\times X$ is trivial, too.
\end{itemize}
The space ${\frak R}$ is the correct parameter space and parameterizes
a universal family $(\E_{\frak R,v}, v\in V; f_{\frak R,a}, a\in A; \EL_{\frak R},
\eps_{\frak R})$. It comes with an action of $\mathop{\sf X}_{v\in V}\GL(W_v)$,
and the universal family is linearized w.r.t.\ that group action.
The parameter space is also projective over $\mathop{\sf X}_{v\in V}\frak Q^0_v$.
\par
The associated family of $V$-split sheaves
$(\E_{\frak R,v}, v\in V; \EL_{\frak R},
\tau_{\frak R})$ of type $\ul{P}$ with a decoration of type $(s,b+1,1,m)$ defines a
$(\mathop{\sf X}_{v\in V}\GL(W_v))$-equivariant morphism
$$
I\colon {\frak R}\lra {\frak M}
$$
over the base scheme $\mathop{\sf X}_{v\in V}\frak Q^0_v$.
Since ${\frak R}$ is proper over $\mathop{\sf X}_{v\in V}\frak Q^0_v$,
the morphism $I$ is automatically proper (\cite{Ha}, II, Cor.~4.8~(e)).
It is also injective. To see this, let $r\in {\frak R}$ be a point which
corresponds to the representation $(\E_v, v\in V; f_a, a\in A; \eps)$
of $Q$ of type $(\ul{P},\ul{\G})$. For any $a\in A$, the surjection
$\Oh_X(-m)^{\oplus b}\otimes W_{t(a)}\otimes\Oh_X(-l)\lra \G_{a}\otimes
\E_{t(a)}$ yields an injective homomorphism
\begin{eqnarray*}
\Hom\bigl(\G_{a}\otimes
\E_{t(a)},\E_{h(a)}\bigr) &\subset&\Hom\bigl(\Oh_X(-m)^{\oplus b}\otimes W_{t(a)}\otimes\Oh_X(-l),
\E_{h(a)}\bigr) \q=
\\
&=&
H^0\bigl({W_{t(a)}^\vee}^{\oplus b}\otimes\E_{h(a)}(l+m)\bigr).
\end{eqnarray*}
Since $\E_{h(a)}$ is torsion free, the restriction map
$$
H^0\bigl(X,{W_{t(a)}^\vee}^{\oplus b}\otimes\E_{h(a)}(l+m)\bigr)
\q\lra\q
H^0\bigl(U,{W_{t(a)}^\vee}^{\oplus b}\otimes\E_{h(a)}(l+m)_{|U}\bigr)
$$
is injective for any open subset $U=X\setminus Z$ with ${\rm codim}_X(Z)\ge 2$.
If ${\frak U}\subset {\frak R}\times X$ is the maximal open subset
where all the $\E_{\frak R,v}$, $v\in V$, are locally free, then
$U:={\frak U}\cap (\{r\}\times X)$ is the maximal open subset where all
the $\E_v$, $v\in V$, are locally free (\cite{HL}, Lemma~2.1.7). In particular, the
complement of $U$ in $X$ has codimension at least two. Since $\tau_{\frak R|\{q\}\times X}$
determines all the $f_a$, $a\in A$, and $\eps$ over $U$, we are done.
Because $I$ is injective and proper and, thus, finite, Theorem~\ref{main2}
follows immediately from Theorem~\ref{main1} and its proof.\qed
\subsection{Proof of Theorem~\ref{stab}}
The proof of Theorem~\ref{stab} is basically a formal adaptation of the corresponding result
for Hitchin pairs.
If, in the following, a representation $(\E_v,v\in V; f_a, a\in A;\eps)$ is given,
$\tau$ will always stand for the associated decoration.
We first observe
\begin{Lem}
\label{bound}
Suppose we are given $\ul{\sigma}$, $\ul{\eta}$, and $\delta$, as well as
a $(\ul{\sigma},\ul{\eta},\delta)$-semistable representation
$(\E_v,v\in V;f_a,a\in A)$ of type
$(\ul{P},\ul{\G})$. Then, it satisfies Condition {\rm 2.a.}\ of Theorem~{\rm\ref{stab}}.
\end{Lem}
\noindent\it Proof\rm.
Let $(\F_v,v\in V)$ be a non-trivial, proper
subrepresentation of $(\E_v,v\in V;f_a,a\in A)$.
Set
$$
\E_\bullet\colon\qquad 0\subsetneq (\F_v,v\in V)\subsetneq (\E_v,v\in V).
$$
Then, one verifies $\mu(\E_\bullet, (1);\tau)\le 0$, from which the assertion follows.
\qed
\par
\begin{Prop}
Fix $\ul{\sigma}$ and $\ul{\chi}$.
Then, the set of torsion free sheaves occurring
in representations of type $(\ul{P},\ul{G})$ which satisfy Condition {\rm 2.a.} of
Theorem~{\rm\ref{stab}} is bounded.
\end{Prop}
\noindent\it Proof\rm.  We fix surjections
$\nu_a\colon \Oh_X(-m)^{\oplus b}\lra\G_a$, $a\in A$.
We may now adapt Nitsure's argument \cite{Nit}, Proposition~3.2.
Let
$$
0=(\F_0^v,v\in V)\subsetneq (\F_1^v,v\in V)\subsetneq\cdots\subsetneq (\F_s^v,v\in V)
\subsetneq (\F_{s+1}^v,v\in V)=(\E_v,v\in V),
$$
be the Harder-Narasimhan filtration of $(\E_v,v\in V)$ defined w.r.t.\ the parameters
$\ul{\sigma}$ and $\ul{\chi}$.
It will suffice to bound $\mu_{\ul{\sigma},\ul{\chi}}(\F_1^v,v\in V)$.
Define
$$
D:= \deg\Oh_X(m)+\max_{v\in V}\left\{\frac{\ol{\chi}_{v}}{\sigma_v}\right\}
-
\min_{v\in V}\left\{\frac{\ol{\chi}_{v}}{\sigma_v}\right\}.
$$
We claim that
$$
\mu_{\ul{\sigma},\ul{\chi}}\bigl(\F_1^v,v\in V\bigr)
\le\min
\left\{\,
\mu_{\ul{\sigma},\ul{\chi}}\bigl(\E_v,v\in V\bigr),
\
\mu_{\ul{\sigma},\ul{\chi}}\bigl(\E_v,v\in V\bigr)
+\frac{\bigl((\sum_{v\in V}r_v)-1\bigr)^2}{\sum_{v\in V}r_v}\cdot D \,\right\}.
$$
We can view the collection $f_a$, $a\in A$, together with the zero
homomorphisms $\E_{v_1}\lra\E_{v_2}\otimes\Oh_X(m)^{\oplus b}$ for
$(v_1,v_2)\not\in A$ as a homomorphism
$$
f\colon\bigoplus_{v\in V} \E_v\lra \Bigl(\bigoplus_{v\in V}\E_v\Bigr)\otimes\Oh_X(m)^{\oplus b}.
$$
For any $V$-split subsheaf $(\F_v,v\in V)$, the condition of being a
subrepresentation is, thus, equivalent to
$f(\bigoplus_{v\in V}\F_v)\subseteq (\bigoplus_{v\in V}\F_v)\otimes\Oh_X(m)^{\oplus b}$.
We simply say that $\bigoplus_{v\in V}\F_v$ is \it $f$-invariant\rm.
If the condition
$$
\mu_{\ul{\sigma},\ul{\chi}}\bigl(\F_1^v,v\in V\bigr)
\le \mu_{\ul{\sigma},\ul{\chi}}\bigl(\E_v,v\in V\bigr)
$$
is violated, then, by definition, none of the sheaves $\F_j^{\rm
total}:=\bigoplus_{v\in V}\F_j^v$ can be $f$-invariant,
$j=1,...,s$, i.e., the homomorphisms $\phi_j\colon  \F_j^{\rm
total}\lra \bigl(\E^{\rm total}/\F_j^{\rm
total}\bigr)\otimes\Oh_X(m)^{\oplus b}$ are non-trivial, $\E^{\rm
total}:=\bigoplus_{v\in V}\E_v$, $j=1,...,s$. For any
$j=1,...,s$, there exist $\iota\le j-1$ and $\kappa\ge j+1$, such
that $\phi_j$ induces a non trivial homomorphism
$\ol{\phi}_j\colon  \F_{\iota+1}^{\rm total}/\F^{\rm
total}_{\iota} \lra \bigl(\F_{\kappa}^{\rm total}/\F^{\rm
total}_{\kappa-1}\bigr)\otimes\Oh_X(m)^{\oplus b}$. Now,
Proposition~\ref{homtriv} and \ref{HS} imply
\begin{eqnarray*}
\mu_{\ul{\sigma},\ul{\chi}}
\bigl(\F^v_j/\F^v_{j-1},v\in V\bigr)
&\le&
\mu_{\ul{\sigma},\ul{\chi}}
\bigl(\F^v_{\iota+1}/\F^v_\iota,v\in V\bigr)
\le
\mu_{\ul{\sigma},\ul{\chi}}
\bigl(\F^v_\kappa/\F^v_{\kappa-1},v\in V\bigr)+D
\\
&\le&
\mu_{\ul{\sigma},\ul{\chi}}
\bigl(\F^v_{j+1}/\F^v_j,v\in V\bigr)+D.
\end{eqnarray*}
Thus,
\begin{eqnarray*}
\mu_{\ul{\sigma},\ul{\chi}}
\bigl(\F^v_1,v\in V\bigr)
&\le& \mu_{\ul{\sigma},\ul{\chi}}
\bigl(\E_v/\F^v_{s},v\in V\bigr)+s\cdot D
\\
&\le & \mu_{\ul{\sigma},\ul{\chi}}
\bigl(\E_v/\F^v_{s},v\in V\bigr)+\Bigl(\bigl(\sum_{v\in V} r_v\bigr)-1\Bigr)\cdot D.
\end{eqnarray*}
Finally, one finds that
$$
\mu_{\ul{\sigma},\ul{\chi}}
\bigl(\E_v/\F^v_{s},v\in V\bigr)\le
\frac{\deg_{\ul{\sigma},\ul{\chi}}(\E_v,v\in V)
-\mu_{\ul{\sigma},\ul{\chi}}
\bigl(\F^v_1,v\in V\bigr)}
{(\sum_{v\in V} r_v)-1}
$$
from which the assertion follows.
\qed
\par
We first show that 1.\ implies 2. We have already checked that 2.a.\ holds.
Before we check Condition 2.b., we review the Linear Algebra setting.
The space $\bigoplus_{a\in A}\Hom\bigl(\C^{r_{t(a)}},\allowbreak \C^{r_{h(a)}}\bigr)^{\oplus b}$
can be $\bigl(\mathop{\sf X}_{v\in V}\GL(\C^{r_v})\bigr)$-equivariantly embedded
into
$$
{\Bbb E}:=\End(M)^{\oplus b},\qquad
M:=\bigoplus_{v\in V}\C^{\sigma_v\cdot r_v}.
$$
Let $\Bigl[f_a\colon  \C^{r_{t(a)}}\lra\bigl(\C^{r_{h(a)}}\bigr)^{\oplus b}, a\in A\Bigr]$
in $\P\Bigl(\bigoplus_{a\in A}\Hom\bigl(\C^{r_{t(a)}},\C^{r_{h(a)}}\bigr)^{\oplus b}\Bigr)$
be an element which is unstable w.r.t.\ $\widetilde{G}$-action. Let
$[f]\in \P({\Bbb E})$ be the associated element which is equally unstable.
As explained before, a one parameter subgroup $\la\colon  \C^*\lra \widetilde{G}$
yields a weighted flag
$(M_\bullet,\ul{\alpha})$ with
$$
M_\bullet\colon\qquad 0\subsetneq U_1\subsetneq \cdots\subsetneq U_s
\subsetneq \C^{\sum_{v\in V}\sigma_v\cdot r_v}.
$$
Here, $U_j=\bigoplus_{v\in V}U_j^{v,\oplus\sigma_v}$ for
suitable subspaces $U^v_j\subset \C^{r_v}$, $j=1,...,s$.
Then, $\mu(\la,[f])<0$ will occur if and only if
$$
f\bigl(U_j\bigr) \subset  U_{j-1}^{\oplus b},\q j=1,...,s+1,
$$
i.e.,
$$
f\Bigl(\bigoplus_{v\in V} U_j^{v,\oplus \sigma_v}\Bigr)
\subset
\Bigl(\bigoplus_{v\in V}
U_{j-1}^{v,\oplus \sigma_v}\Bigr)^{\oplus b},\qquad j=1,...,s+1,
$$
or equivalently
$$
U_j^{v,\oplus\sigma_v} \subset  \ker\biggl(M
\stackrel{f}{\lra} \bigoplus_{v\in V}
\Bigl(\C^{\sigma_v\cdot r_v}/U_{j-1}^{v,\oplus\sigma_v}\Bigr)^{\oplus b}
\biggr),
\q\hbox{for all } v\in V.
$$
Conversely, we may define
$$
Y^v_1:=\ker\biggl(\C^{r_v}\stackrel{\rm diag}{\hookrightarrow}\C^{\sigma_v\cdot r_v}
\stackrel{f}{\lra}M^{\oplus b}\biggr),\q v\in V,
$$
and
$$
Y^v_j:=\ker\biggl(\C^{r_v}\stackrel{\rm diag}{\hookrightarrow}\C^{\sigma_v\cdot r_v}
\stackrel{f}{\lra}
\bigoplus_{v\in V}
\Bigl(\C^{\sigma_v\cdot r_v}/Y_{j-1}^{v,\oplus\sigma_v}\Bigr)^{\oplus b}\biggr),\q v\in V.
$$
By our previous observations, this process will stop after at most $s$ steps, i.e., we get
a flag
$$
M^\p_\bullet\colon\qquad
0\subsetneq \bigoplus_{v\in V} Y^{v,\oplus\sigma_v}_1
\subsetneq \cdots\subsetneq \bigoplus_{v\in V} Y^{v,\oplus\sigma_v}_{s^\p}
\subsetneq \C^{\sum_{v\in V}\sigma_v\cdot r_v},
$$
and $(M^\p_\bullet, (1,...,1))$ comes from a suitable one parameter subgroup
$\la\colon \C^*\lra \widetilde{G}$ with $\mu(\la,[f])=-\sum_{v\in V}\sigma_v\cdot r_v$.
\par
The latter construction can be extended to the setting of sheaves, i.e.,
given a representation $(\E_v,v\in V;f_a, a\in A;\eps=0)$ as in 1.\ for which
Condition 2.b.\ fails, we define
$$
\F^v_1:=\ker\biggl(\E_v\stackrel{\rm diag}{\hookrightarrow}\E_v^{\oplus\sigma_v}
\stackrel{f}{\lra}\bigoplus_{v\in V}
\E_v^{\oplus\sigma_v}\otimes\Oh_X(m)^{\oplus b}\biggr),\q v\in V,
$$
and
$$
\F^v_j:=\ker\biggl(\E_v\stackrel{\rm diag}{\hookrightarrow}\E_v^{\oplus\sigma_v}
\stackrel{f}{\lra}
\bigoplus_{v\in V}
\Bigl(\E_v^{\oplus\sigma_v}/\F_{j-1}^{v,\oplus\sigma_v}\Bigr)\otimes\Oh_X(m)^{\oplus b}\biggr),\q v\in V,
j>1.
$$
Then, we find the weighted filtration $(\E^\bullet,(1,...,1))$ with
$$
\E_\bullet\colon\qquad 0\subsetneq (\F^v_1,v\in V)\subsetneq
\cdots\subsetneq (\F^v_{s^\p},v\in V) \subsetneq (\E_v, v\in V)
$$
and
$$
\mu\bigl(\E^\bullet,(1,...,1);\tau\bigr)= -\sum_{v\in V}\sigma_v\cdot r_v.
$$
By the Boundedness Result \ref{bound}, it is clear that the
sheaves of the form $\F_j^v$ as just defined live in bounded
families, too. In particular, there is a constant $C>0$, such that
$$
\deg_{\ul{\sigma},\ul{\chi}}(\E_v,v\in V)\cdot{\rk}_{\ul{\sigma}}(\F_j^v,v\in V)
-
\deg_{\ul{\sigma},\ul{\chi}}(\F_j^v,v\in V)\cdot{\rk}_{\ul{\sigma}}(\E_v,v\in V)
< C
$$
for any filtration as above.
But then, with $\ol{\delta}>0$, the coefficient of $x^{\dim X-1}$ in $\delta$,
the condition of $(\ul{\sigma},\ul{\eta}/n,n\cdot\delta)$-semistability requires
\begin{eqnarray*}
0&\preceq& M_{\ul{\sigma},\ul{\chi}}\bigl(\E_\bullet,(1,...,1)\bigr)
+n\cdot\delta\cdot\mu\bigl(\E_\bullet,(1,...,1);\tau\bigr)
\\
 &\preceq& \Bigl(\bigl((\sum_{v\in V} r_v)-1\bigr)
 \cdot C-n\cdot\ol{\delta}\cdot\sum_{v\in V}(\sigma_v\cdot
 r_v)\Bigr)\cdot x^{\dim X-1},
\end{eqnarray*}
but for large $n$, this is impossible.
\par
The converse is an easy adaptation of the argument given in \cite{Sch}, Example 3.6,
and is left as an exercise to the reader. \qed
\subsection{Proof of Theorem~\ref{main3}}
The points i) and ii) are just a reformulation of Theorem~\ref{main2}.
Point iii) is proved by standard arguments and will be omitted here.
Finally, the constructions carried out in Section~\ref{large} show that the universal
family on the parameter space ${\frak R}$ defines a morphism
${\frak R}\lra {\Bbb H}(Q,\ul{P},\ul{\G})$. This morphism is invariant under
the $\widetilde{G}$-action and, thus, descends to the moduli space
${\cal R}(Q)^{(\ul{\sigma},\ul{\chi})-\rm ss}_{\ul{P}/\ul{\G}}$.\qed
\section{Appendix: Two auxiliary results}
\subsection{Restrictions of families of locally free sheaves to open subsets}
Let $X$ be a smooth projective manifold and $S$ a noetherian scheme. Let $\iota\colon
{U}\subset S\times X$
be an open subset, such that
$$
{\rm codim}\bigl(X\setminus ({U}\cap \{s\}\times X), X\bigr)\ge 2,\q
\hbox{for all $s\in S$}.
$$
\begin{Prop}
\label{Mar}
In the above situation, the natural homomorphism
$\Oh_{S\times X}\lra\iota_*\Oh_{U}$ is an isomorphism. In particular, for any
locally free sheaf ${\cal V}$ on $S\times X$, we have
$$
\cal V= \iota_*\iota^*\cal V.
$$
\end{Prop}
\noindent
\it Proof\rm.
We refer to \cite{Ma2}, p.~111f.\qed
\subsection{Zero loci of sheaf homomorphisms}
The following result may be found in \cite{GS}, Lemma~3.1.
\begin{Prop}
\label{GS}
Let $S$ be a noetherian scheme,
$\cal A^1_S$ and $\cal A_S^2$ coherent sheaves on
$S\times X$, and $\phi_S\colon \cal A_S^1\lra \cal A_S^2$
a homomorphism. Assume that $\cal A_S^2$ is $S$-flat.
Then, there is a closed subscheme ${\frak Y}\subset S$ the closed points of which
are
those $s\in S$ for which $\phi_{S|\{s\}\times X}\equiv 0$.
More precisely, it has the property that any morphism $f\colon T\lra S$ factors
through ${\frak Y}$, if and only if $(f\times \id_X)^*\phi_S$ is the zero
homomorphism.
\end{Prop}
\section{Appendix: A concluding remark} 
It was pointed out to me by Bal\'azs Szendr\"oi that the case in which $X$ is a point is formally
not covered by our formalism.
It would be formally included, if one allowed $\delta$ to have degree
$\dim(X)$ (or higher). Then, for $X$ a point, $\delta=1$, $\eta_v$, $v\in V$, with $\sum_{v\in V}\eta_v r_v=0$, and
a representation $(f_a, a\in A)\in \bigoplus_{a\in A}\Hom\bigl(\C^{r_{t(a)}},\C^{r_{h(a)}}\bigr)$,
the condition of (semi)stability would
read
\begin{eqnarray*}
&&
\frac{\sum_{v\in V}\sigma_v\dim(W_v)-\sum_{v\in V}\eta_v \dim(W_v)}
{\sum_{v\in V}\sigma_v\dim(W_v)}
=1-\frac{\sum_{v\in V}\eta_v \dim(W_v)}
{\sum_{v\in V}\sigma_v\dim(W_v)}
\\
&(\le)& \frac{\sum_{v\in V}\sigma_vr_v-\sum_{v\in V}\eta_v r_v}
{\sum_{v\in V}\sigma_vr_v}
=1,
\end{eqnarray*}
i.e.,
$$
-\sum_{v\in V}\eta_v \dim(W_v)(\le)0
$$
for any non-trivial subrepresentation $(W_v,v\in V)$. This is precisely King's definition
w.r.t.\ the character
$$
\bigl(U_v,v\in V;f_a, a\in A)\lma -\sum_{v\in V} \eta_v\cdot \dim(U_v).
$$
Note that for $\dim(X)\ge 0$, a positive polynomial $\delta$ of degree
$\dim(X)$, and $\eta_v$, $v\in V$, with $\sum_{v\in V}\eta_vr_v=0$, the condition
of (semi)stability can be restated as follows:
\begin{itemize}
\item For any non-trivial subrepresentation $(\F_v,v\in V)$ one has
$$
-\sum_{v\in V} \eta_v\rk(\F_v)(\le)0,
$$
and,
\item if ``$=$" occurs, then
$$
\frac{\sum_{v\in V}\sigma_v P(\F_v)}{\sum_{v\in V} \sigma_v\rk(\F_v)}\q (\preceq)\q
\frac{\sum_{v\in V}\sigma_v P(\E_v)}{\sum_{v\in V} \sigma_v\rk(\E_v)}.
$$
\end{itemize}
Let us call representations which satisfy this condition \it asymptotically (semi)stable\rm.
The moduli spaces for asymptotically (semi)stable objects might be obtained as follows:
\begin{itemize}
\item Fix the data $\sigma_v$, $\eta_v$, $v\in V$.
\item There is a positive polynomial $\delta_0$ (depending on the Hilbert polynomials
$P_v$, $\sigma_v$, and $\eta_v$, $v\in V$) of degree $\dim(X)-1$, such that, for any
$\delta\succ\delta_0$, a representation $(\E_v, v\in V;f_a, a\in A)$ will be (semi)stable
w.r.t.\ the parameters $\sigma_v$, $\eta_v$, $v\in V$, and $\delta$, if and only if it is
asymptotically (semi)stable.
\end{itemize}
The techniques to prove this should be adapted from my recent paper \cite{SchGlobal}.

\end{document}